\documentclass[11pt]{amsart}

\usepackage{kbordermatrix}
\usepackage{graphicx}
\usepackage{amssymb}
\usepackage{times}
\usepackage[top=25mm, bottom=25mm, left=25mm, right=25mm]{geometry}
\usepackage{color}
\usepackage{bm}
\newenvironment{widequote}{\begin{list}{}
      {\setlength{\rightmargin}{0mm}\setlength{\leftmargin}{5mm}}
      \item[]}{\end{list}}

\newcommand{\pgnotate}[1]{}

\newcommand{\vct}[1]{\bm{\mathsf{#1}}}
\newcommand{\mtx}[1]{\bm{\mathsf{#1}}}

\newtheorem{thm}{Theorem}
\newtheorem{lemma}[thm]{Lemma}

\theoremstyle{definition}
\newtheorem{remark}{Remark}

\usepackage{fancyhdr}
\begin{document}

\begin{center}
\textbf{\large{Randomized methods for matrix computations}}

\vspace{5mm}

Per-Gunnar Martinsson

\vspace{3mm}

Institute for Computational Sciences and Engineering\\
The University of Texas at Austin

\vspace{3mm}

January 31, 2019

\vspace{20mm}

\begin{minipage}{130mm}
\textbf{Abstract:}
The purpose of this text is to provide an accessible introduction   
to a set of recently developed algorithms for factorizing matrices.
These new algorithms attain high practical speed by reducing the
dimensionality of intermediate computations using randomized projections.
The algorithms are particularly powerful for computing low-rank approximations
to very large matrices, but they can also be used to accelerate algorithms
for computing full factorizations of matrices. A key competitive advantage
of the algorithms described is that they require less communication than
traditional deterministic methods.
\end{minipage}
\end{center}

\clearpage

\tableofcontents

\clearpage

\section{Introduction}

\subsection{Scope and objectives}

The objective of this text is to describe a set of randomized
methods for efficiently computing a low rank approximation to a given
matrix. In other words, given an $m\times n$ matrix $\mtx{A}$, we seek
to compute factors $\mtx{E}$ and $\mtx{F}$ such that
\begin{equation}
\label{eq:generalEF}
\begin{array}{ccccccccccccccc}
\mtx{A} &\approx& \mtx{E} & \mtx{F},\\
m\times n   && m\times k & k\times n
\end{array}
\end{equation}
where the rank $k$ of the approximation is a number we assume to be much
smaller than either $m$ or $n$. In some situations, the rank $k$ is given
to us in advance, while in others, it is part of the problem to determine
a rank such that the approximation satisfies a bound of the type
$$
\|\mtx{A} - \mtx{E}\mtx{F}\| \leq \varepsilon
$$
where $\varepsilon$ is a given tolerance, and $\|\cdot\|$ is some specified
matrix norm (in this text, we will discuss only the spectral and the
Frobenius norms).

An approximation of the form (\ref{eq:generalEF}) is useful for
storing the matrix $\mtx{A}$ more frugally (we can store $\mtx{E}$ and
$\mtx{F}$ using $k(m+n)$ numbers, as opposed to $mn$ numbers for storing
$\mtx{A}$), for efficiently computing a matrix vector product $\vct{z} = \mtx{A}\vct{x}$
(via $\vct{y} = \mtx{F}\vct{x}$ and $\vct{z} = \mtx{E}\vct{y}$), for
data interpretation, and much more.
Low-rank approximation problems of this type form a cornerstone of data analysis
and scientific computing, and arise in a broad range of applications,
including principal component analysis (PCA) in computational statistics,
spectral methods for clustering high-dimensional data and finding structure
in graphs, image and video compression, model reduction in physical modeling,
and many more.

In performing low-rank approximation, one is typically interested in specific
factorizations where the factors $\mtx{E}$ and $\mtx{F}$ satisfy additional
constraints. When $\mtx{A}$ is a symmetric $n\times n$ matrix, one is commonly interested in
finding an approximate rank-$k$ eigenvalue decomposition (EVD), which takes the form
\begin{equation}
\begin{array}{ccccccccccccccc}
\mtx{A} &\approx& \mtx{U} & \mtx{D} & \mtx{U}^{*},\\
n\times n   && n\times k & k\times k & k\times n
\end{array}
\end{equation}
where the columns of $\mtx{U}$ form an orthonormal set, and where $\mtx{D}$ is
diagonal. For a general $m\times n$ matrix $\mtx{A}$, we would typically be
interested in an approximate rank-$k$ singular value decomposition (SVD), which
takes the form
\begin{equation}
\begin{array}{ccccccccccccccc}
\mtx{A} &\approx& \mtx{U} & \mtx{D} & \mtx{V}^{*},\\
m\times n   && m\times k & k\times k & k\times n
\end{array}
\end{equation}
where $\mtx{U}$ and $\mtx{V}$ have orthonormal columns, and $\mtx{D}$ is diagonal.
In this text, we will discuss both the EVD and the SVD in depth. We will also
describe factorizations such as the \textit{interpolative decomposition (ID)}, and
the \textit{CUR decomposition} which are highly useful for data interpretation, and
for certain applications in scientific computing. In these factorizations, we seek
to determine a subset of the columns (rows) of $\mtx{A}$ itself that together form a good
(approximate) basis for the column (row) space.

While most of the text is aimed at computing low rank factorizations where the
target rank $k$ is much smaller than the dimensions of the matrix $m$ and $n$, we
will in the last couple of sections of the text also discuss how randomization
can be used to speed up factorization of \textit{full} matrices, such as a full
column pivoted QR factorization, or various relaxations of the SVD that are useful
for solving least-squares problems, etc.

\subsection{The key ideas of randomized low-rank approximation}

To quickly introduce the central ideas, let us describe a
simple prototypical randomized algorithm: Let $\mtx{A}$ be a matrix of size $m\times n$
that is approximately of low rank. In other words, we assume that for some integer
$k < \min(m,n)$, there exists an approximate
low rank factorization of the form (\ref{eq:generalEF}).
Then a natural question is how do you in a computationally efficient manner
construct the factors $\mtx{E}$ and $\mtx{F}$? In \cite{2006_martinsson_random1_orig},
it was observed that random matrix theory provides a simple solution:
Draw a \textit{Gaussian random matrix} $\mtx{G}$ of size $n\times k$,
form the \textit{sampling matrix}
$$
\mtx{E} = \mtx{A}\mtx{G},
$$
and then compute the factor $\mtx{F}$ via
$$
\mtx{F} = \mtx{E}^{\dagger}\mtx{A},
$$
where $\mtx{E}^{\dagger}$ is the Moore-Penrose pseudo-inverse of $\mtx{A}$, cf.~Section \ref{sec:pinv}.
(Then $\mtx{E}\mtx{F} = \mtx{E}\mtx{E}^{\dagger}\mtx{A}$, where
$\mtx{E}\mtx{E}^{\dagger}$ is the orthogonal projection onto the linear
space spanned by the $k$ columns in $\mtx{E}$.)
Then in many important situations, the approximation
\begin{equation}
\label{eq:def_Ak}
\begin{array}{ccccccccccccccc}
\mtx{A} &\approx& \mtx{E} & \bigl(\mtx{E}^{\dagger}\mtx{A}\bigr), \\
m\times n   && m\times k & k\times n
\end{array}
\end{equation}
is close to optimal. With this observation as a starting point, we will
construct highly efficient algorithms for computing approximate spectral
decompositions of $\mtx{A}$, for solving certain least-squares problems,
for doing principal component analysis of large data sets, etc.


\subsection{Advantages of randomized methods}

The algorithms that result from using randomized sampling techniques are
computationally efficient, and are simple to implement as they rely on
standard building blocks (matrix-matrix multiplication, unpivoted QR
factorization, etc.) that are readily available for most computing environments
(multicore CPU, GPU, distributed memory  machines, etc). As an illustration,
we invite the reader to peek ahead at Figure \ref{fig:RSVD}, which provides
a complete Matlab code for a randomized algorithm that computes an approximate
singular value decomposition of a matrix.
Examples of improvements enabled by these randomized algorithms include:
\begin{itemize}
\item Given an $m\times n$ matrix $\mtx{A}$, the cost of computing a rank-$k$ approximant
using classical methods is $O(mnk)$. Randomized algorithms can attain
complexity $O(mn\log k + k^{2}(m+n))$ \cite{2011_martinsson_randomsurvey}*{Sec.~6.1}, and Section \ref{sec:SRFT}.
\item Algorithms for performing principal component analysis (PCA) of large data sets
have been greatly accelerated, in particular when the data is stored out-of-core \cite{2011_halko_largedata}.
\item Randomized methods tend to require less communication than traditional methods, and
can be efficiently implemented on severely communication constrained environments such as
GPUs \cite{2015_martinsson_blocked} and distributed computing platforms
\cite{2012_halko_dissertation}*{Ch.~4} and
\cite{url_facebook,gittens2016matrix}.
\item Randomized algorithms have enabled the development of \textit{single-pass}
matrix factorization algorithms in which
the matrix is ``streamed'' and never stored, cf.~\cite{2011_martinsson_randomsurvey}*{Sec.~6.3} and Section \ref{sec:singlepass}.
\end{itemize}

\subsection{Very brief literature survey}
Our focus in this text is to describe randomized methods that attain high practical
computational efficiency. In particular, we use randomization mostly as a tool for
minimizing \textit{communication}, rather than minimizing the flop count (although we
do sometimes improve asymptotic flop counts as well). The methods described were
first published in \cite{2006_martinsson_random1_orig} (which was inspired by
\cite{2004_kannan_vempala}, and later led to
\cite{2007_martinsson_PNAS,2011_martinsson_random1};
see also \cite{2006_sarlos_improved}). Our presentation largely follows that in the
2011 survey \cite{2011_martinsson_randomsurvey}, but with a focus more on practical
usage, rather than theoretical analysis. We have also included material from more recent
work, including
\cite{2014_martinsson_CUR} on factorizations that allow for better data interpretation,
\cite{2015_martinsson_blocked} on blocking and adaptive error estimation, and
\cite{2017_martinsson_UTV,2015_blockQR_SISC} on full factorizations.

The idea of using randomization to improve algorithms for low-rank approximation of
matrices has been extensively investigated within the theoretical computer science
community, with early work including \cite{2006_drineas_kannan_mahoney,2004_kannan_vempala,2006_sarlos_improved}.
The focus of these texts has been to develop algorithms with optimal or close to optimal theoretical performance guarantees in terms of asymptotic flop counts and error bounds.
The surveys \cite{2011_mahoney_survey,2018_PCMI_mahoney_drineas} and \cite{2014_woodruff_sketching} provide
introductions to this literature.

A version of this tutorial has appeared in \cite{2018_PCMI_martinsson}.

\section{Notation}

\subsection{Notation}
\label{sec:notation}
Throughout the text, we measure vectors in $\mathbb{R}^{n}$ using their
Euclidean norm, $\|\vct{v}\| = \sqrt{\sum_{j=1}^{n}|\vct{v}(i)|^{2}}$. We measure matrices using the spectral and the Frobenius norms,
defined by
$$
\|\mtx{A}\| = \sup_{\|\vct{x}\|=1}\|\mtx{A}\vct{x}\|,
\qquad\mbox{and}\qquad
\|\mtx{A}\|_{\rm Fro} = \left(\sum_{i,j}|\mtx{A}(i,j)|^{2}\right)^{1/2},
$$
respectively.
We use the notation of Golub and Van Loan \cite{golub} to specify
submatrices. In other words, if $\mtx{B}$ is an $m\times n$ matrix with
entries $\mtx{B}(i,j)$, and $I = [i_{1},\,i_{2},\,\dots,\,i_{k}]$ and $J =
[j_{1},\,j_{2},\,\dots,\,j_{\ell}]$ are two index vectors, then
$\mtx{B}(I,J)$ denotes the $k\times \ell$ matrix
$$
\mtx{B}(I,J) = \left[\begin{array}{cccc}
B(i_{1},j_{1}) & B(i_{1},j_{2}) & \cdots & B(i_{1},j_{\ell}) \\
B(i_{2},j_{1}) & B(i_{2},j_{2}) & \cdots & B(i_{2},j_{\ell}) \\
\vdots         & \vdots         &        & \vdots            \\
B(i_{k},j_{1}) & B(i_{k},j_{2}) & \cdots & B(i_{k},j_{\ell})
\end{array}\right].
$$
We let $\mtx{B}(I,:)$ denote the matrix $\mtx{B}(I,[1,\,2,\,\dots,\,n])$,
and define $\mtx{B}(:,J)$ analogously.

The transpose of $\mtx{B}$ is denoted $\mtx{B}^{*}$, and we say that a matrix $\mtx{U}$
is \textit{orthonormal (ON)} if its  columns form an orthonormal set, so that $\mtx{U}^{*}\mtx{U} = \mtx{I}$.

\subsection{The singular value decomposition (SVD)}
\label{sec:SVD}
The SVD was introduced briefly in the introduction. Here we define it again, with some
more detail added. Let $\mtx{A}$ denote an $m\times n$ matrix, and set $r = \min(m,n)$.
Then $\mtx{A}$ admits a factorization
\begin{equation}
\label{eq:svd}
\begin{array}{ccccccccccccccccccccc}
\mtx{A} &=& \mtx{U} & \mtx{D} & \mtx{V}^{*},\\
m\times n && m\times r & r \times r & r\times n
\end{array}
\end{equation}
where the matrices $\mtx{U}$ and $\mtx{V}$ are orthonormal, and $\mtx{D}$ is diagonal.
We let $\{\vct{u}_{i}\}_{i=1}^{r}$ and $\{\vct{v}_{i}\}_{i=1}^{r}$ denote the columns of
$\mtx{U}$ and $\mtx{V}$, respectively. These vectors are the left and right singular vectors
of $\mtx{A}$. The diagonal elements $\{\sigma_{j}\}_{j=1}^{r}$ of
$\mtx{D}$ are the singular values of $\mtx{A}$. We order these so that
$\sigma_{1}  \geq \sigma_{2} \geq \cdots \geq \sigma_{r} \geq 0$.
We let $\mtx{A}_{k}$ denote the truncation of the SVD to its first $k$ terms,
$$
\mtx{A}_{k} =
\mtx{U}(:,1:k)\mtx{D}(1:k,1:k)\bigl(\mtx{V}(:,1:k)\bigr)^{*} =
\sum_{j=1}^{k} \sigma_{j}\vct{u}_{j}\vct{v}_{j}^{*}.
$$
It is easily verified that
\begin{equation}
\label{eq:minerrors}
\|\mtx{A} - \mtx{A}_{k}\| = \sigma_{k+1},
\qquad\mbox{and that}\qquad
\|\mtx{A} - \mtx{A}_{k}\|_{\rm Fro} = \left(\sum_{j=k+1}^{\min(m,n)} \sigma_{j}^{2}\right)^{1/2},
\end{equation}
where $\|\mtx{A}\|$ denotes the operator norm of $\mtx{A}$ and
$\|\mtx{A}\|_{\rm Fro}$ denotes the Frobenius norm of $\mtx{A}$. Moreover,
the Eckart-Young theorem \cite{1936_eckart_young} states that these errors are the smallest possible
errors that can be incurred when approximating $\mtx{A}$ by a matrix of rank $k$.


\subsection{Orthonormalization}
\label{sec:orth}
Given an $m\times \ell$ matrix $\mtx{X}$, with $m \geq \ell$, we introduce the function
$$
\mtx{Q} = \texttt{orth}(\mtx{X})
$$
to denote orthonormalization of the columns of $\mtx{X}$. In other words, $\mtx{Q}$ will be
an $m\times \ell$ orthonormal matrix whose columns form a basis for the column space of $\mtx{X}$.
In practice, this step is typically achieved most efficiently by a call to a packaged QR
factorization; e.g., in Matlab, we would write $[\mtx{Q},\sim] = \texttt{qr}(\mtx{X},0)$.
However, all calls to \texttt{orth} in this manuscript can be implemented \textit{without pivoting},
which makes efficient implementation much easier.

\subsection{The Moore-Penrose pseudoinverse}
\label{sec:pinv}
The Moore-Penrose pseudoinverse is a generalization of the concept of an inverse
for a non-singular square matrix. To define it, let $\mtx{A}$ be a given $m\times n$
matrix. Let $k$ denote its actual rank, so that its singular value decomposition (SVD)
takes the form
$$
\mtx{A} =
\sum_{j=1}^{k}\sigma_{j}\vct{u}_{j}\vct{v}_{j}^{*} =
\mtx{U}_{k}\mtx{D}_{k}\mtx{V}_{k}^{*},
$$
where $\sigma_{1} \geq \sigma_{2} \geq \sigma_{k} > 0$.
Then the pseudoinverse of $\mtx{A}$ is the $n\times m$ matrix defined via
$$
\mtx{A}^{\dagger} =
\sum_{j=1}^{k}\frac{1}{\sigma_{j}}\vct{v}_{j}\vct{u}_{j}^{*} =
\mtx{V}_{k}\mtx{D}_{k}^{-1}\mtx{U}_{k}^{*}.
$$
For any matrix $\mtx{A}$, the matrices
$$
\mtx{A}^{\dagger}\mtx{A} = \mtx{V}_{k}\mtx{V}_{k}^{*},
\qquad\mbox{and}\qquad
\mtx{A}\mtx{A}^{\dagger} = \mtx{U}_{k}\mtx{U}_{k}^{*},
$$
are the orthogonal projections onto the row and column spaces of $\mtx{A}$, respectively.
If $\mtx{A}$ is square and non-singular, then $\mtx{A}^{\dagger} = \mtx{A}^{-1}$.

\section{A two-stage approach}
\label{sec:twostage}

The problem of computing an approximate low-rank factorization to a given
matrix can conveniently be split into two distinct ``stages.'' For concreteness,
we describe the split for the specific task of computing an approximate singular
value decomposition. To be precise, given an $m\times n$ matrix $\mtx{A}$
and a target rank $k$, we seek to compute factors $\mtx{U}$, $\mtx{D}$,
and $\mtx{V}$ such that
\[
\begin{array}{cccccccccccccccccc}
\mtx{A} &\approx& \mtx{U} & \mtx{D} & \mtx{V}^{*}.\\
m\times n && m\times k & k\times k & k\times n
\end{array}
\]
The factors $\mtx{U}$ and $\mtx{V}$ should be orthonormal, and $\mtx{D}$
should be diagonal.
(For now, we assume that the rank $k$ is known in advance, techniques for relaxing
this assumption are described in Section \ref{sec:acc_given}.)
Following \cite{2011_martinsson_randomsurvey}, we split
this task into two computational stages:

\begin{description}

\item[Stage A --- find an approximate range]
Construct an $m\times k$ matrix $\mtx{Q}$ with orthonormal columns
such that $\mtx{A} \approx \mtx{Q}\mtx{Q}^{*}\mtx{A}$.
(In other words, the columns of $\mtx{Q}$ form an approximate
basis for the column space of $\mtx{A}$.) This step
will be executed via a randomized process described in Section \ref{sec:stageA}.

\vspace{2mm}

\item[Stage B --- form a specific factorization]
Given the matrix $\mtx{Q}$ computed in Stage A, form the factors
$\mtx{U}$, $\mtx{D}$, and $\mtx{V}$ using classical deterministic
techniques. For instance, this stage can be executed via the following steps:
\begin{enumerate}
\item Form the $k\times n$ matrix $\mtx{B} = \mtx{Q}^{*}\mtx{A}$.
\item Compute the SVD of the (small) matrix $\mtx{B}$ so that
      $\mtx{B} = \hat{\mtx{U}}\mtx{D}\mtx{V}^{*}$.
\item Form $\mtx{U} = \mtx{Q}\hat{\mtx{U}}$.
\end{enumerate}
\end{description}

The point here is that in a situation where $k \ll \min(m,n)$, the
difficult part of the computation is all in Stage A. Once
that is finished, the post-processing in Stage B is easy, as
all matrices involved have at most $k$ rows or columns.

\begin{remark}
\label{remark8:nomoreerror}
Stage B is exact up to floating point arithmetic so all errors
in the factorization process are incurred at Stage A. To be precise, we have
$$
\mtx{Q}\underbrace{\mtx{Q}^{*}\mtx{A}}_{=\mtx{B}} =
\mtx{Q}\underbrace{\mtx{B}}_{=\hat{\mtx{U}}\mtx{D}\mtx{V}^{*}} =
\underbrace{\mtx{Q}\hat{\mtx{U}}}_{=\mtx{U}}\mtx{D}\mtx{V}^{*} =
\mtx{U}\mtx{D}\mtx{V}^{*}.
$$
In other words, if the factor $\mtx{Q}$ satisfies
$\|\mtx{A} - \mtx{Q}\mtx{Q}^{*}\mtx{A}\| \leq \varepsilon$,
then automatically
\begin{equation}
\label{eq8:error_SVD}
\|\mtx{A} - \mtx{U}\mtx{D}\mtx{V}^{*}\| =
\|\mtx{A} - \mtx{Q}\mtx{Q}^{*}\mtx{A}\| \leq \varepsilon
\end{equation}
unless $\varepsilon$ is close to the machine precision.
\end{remark}

\begin{remark}
A bound of the form (\ref{eq8:error_SVD}) implies that the diagonal
elements $\{\mtx{D}(i,i)\}_{i=1}^{k}$ of $\mtx{D}$ are accurate approximations
to the singular values of $\mtx{A}$ in the sense that $|\sigma_{i} - \mtx{D}(i,i)| \leq \varepsilon$
for $i = 1,\,2,\,\dots,\,k$. However, a bound like (\ref{eq8:error_SVD})
does not provide assurances on the \textit{relative errors} in the singular
values; nor does it provide strong assurances that the columns of $\mtx{U}$
and $\mtx{V}$ are good approximations to the singular vectors of $\mtx{A}$.
\pgnotate{Chase down pointer to literature.}
\end{remark}

\section{A randomized algorithm for ``Stage A'' --- the range finding problem}
\label{sec:stageA}

This section describes a randomized technique for solving
the range finding problem introduced as ``Stage A'' in Section
\ref{sec:twostage}. As a preparation for this discussion, let
us recall that an ``ideal'' basis matrix $\mtx{Q}$
for the range of a given matrix $\mtx{A}$ is the matrix
$\mtx{U}_{k}$ formed by the $k$ leading left singular vectors
of $\mtx{A}$. Letting $\sigma_{j}(\mtx{A})$ denote the $j$'th
singular value of $\mtx{A}$, the Eckart-Young theorem
\cite{1993_stewart_historySVD} states that
\[
\inf\{||\mtx{A} - \mtx{C}||\,\colon\,\mtx{C}\mbox{ has rank }k\} =
||\mtx{A} - \mtx{U}_{k}\mtx{U}_{k}^{*}\mtx{A}|| =
\sigma_{k+1}(\mtx{A}).
\]

Now consider a simplistic randomized method for constructing a spanning
set with $k$ vectors for the range of a matrix $\mtx{A}$:
Draw $k$ random vectors $\{\vct{g}_{j}\}_{j=1}^{k}$ from a Gaussian
distribution, map these to vectors $\vct{y}_{j} = \mtx{A}\vct{g}_{j}$ in the range of $\mtx{A}$,
and then use the resulting set $\{\vct{y}_{j}\}_{j=1}^{k}$
as a basis. Upon orthonormalization via, e.g., Gram-Schmidt, an orthonormal
basis $\{\vct{q}_{j}\}_{j=1}^{k}$ would be obtained. For the special case where
the matrix $\mtx{A}$ has \textit{exact} rank $k$, one can prove that the vectors
$\{\mtx{A}\vct{g}_{j}\}_{j=1}^{k}$ would with probability $1$ be linearly independent,
and the resulting orthonormal (ON) basis $\{\vct{q}_{j}\}_{j=1}^{k}$ would therefore
exactly span the range of $\mtx{A}$. This would in a sense be an ideal algorithm.
The problem is that in practice, there are almost always many non-zero singular values
beyond the first $k$ ones. The left singular vectors associated with these modes all
``pollute'' the sample vectors $\vct{y}_{j} = \mtx{A}\vct{g}_{j}$ and will therefore shift the space
spanned by $\{\vct{y}_{j}\}_{j=1}^{k}$ so that it is no longer aligned with the ideal space
spanned by the $k$ leading singular vectors of $\mtx{A}$. In consequence, the process
described can (and frequently does) produce a poor basis. Luckily, there is
a fix: Simply take a few extra samples. It turns out that if we take, say, $k+10$
samples instead of $k$, then the process
will with probability almost 1 produce a basis that is comparable to the best
possible basis.

\begin{figure}
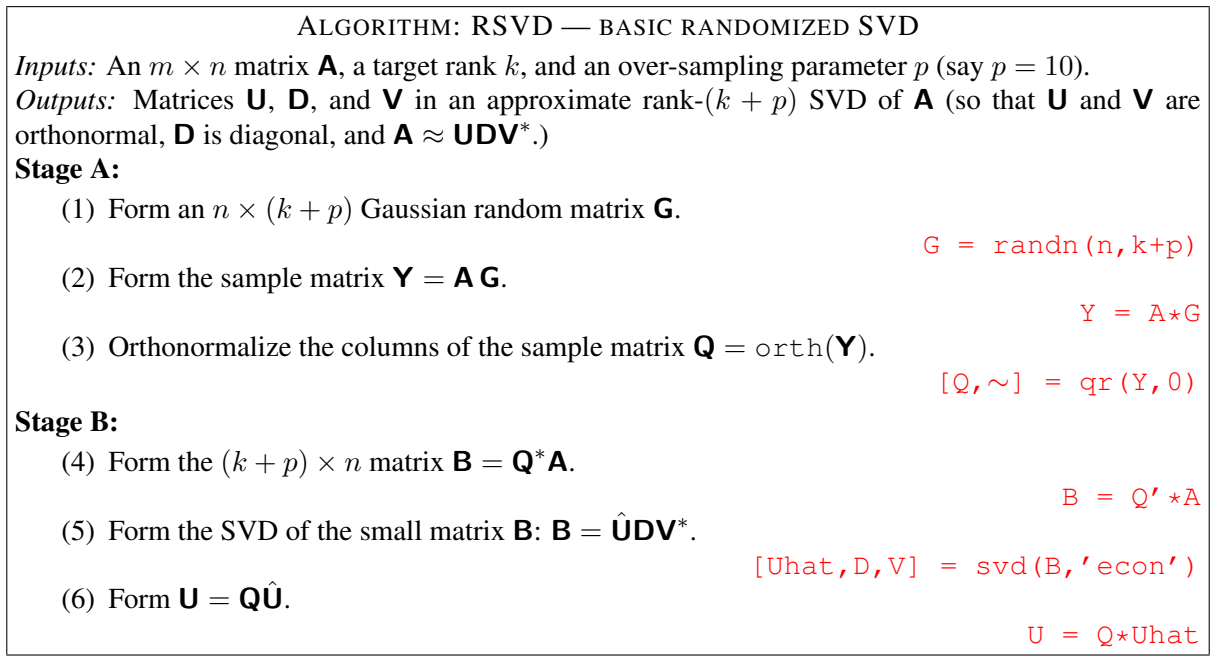

\fbox{\begin{minipage}{0.95\textwidth}
\begin{center}
\textsc{Algorithm: RSVD --- basic randomized SVD}
\end{center}


\textit{Inputs:} An $m\times n$ matrix $\mtx{A}$, a target rank $k$, and an over-sampling parameter $p$ (say $p=10$).


\textit{Outputs:} Matrices $\mtx{U}$, $\mtx{D}$, and $\mtx{V}$ in an approximate rank-$(k+p)$ SVD
of $\mtx{A}$ (so that $\mtx{U}$ and $\mtx{V}$ are orthonormal, $\mtx{D}$ is diagonal, and
$\mtx{A} \approx \mtx{U}\mtx{D}\mtx{V}^{*}$.)


\textbf{Stage A:}
\begin{enumerate}
\item Form an $n\times (k+p)$ Gaussian random matrix $\mtx{G}$.\\
\mbox{} \hfill \texttt{\color{red}G = randn(n,k+p)}
\item Form the sample matrix $\mtx{Y} = \mtx{A}\,\mtx{G}$.\\
\mbox{} \hfill \texttt{\color{red}Y = A*G}
\item Orthonormalize the columns of the sample matrix $\mtx{Q} = \texttt{orth}(\mtx{Y})$.\\
\mbox{} \hfill \texttt{\color{red}[Q,$\sim$] = qr(Y,0)}
\end{enumerate}

\textbf{Stage B:}
\begin{enumerate}
\setcounter{enumi}{3}
\item Form the $(k+p)\times n$ matrix $\mtx{B} = \mtx{Q}^{*}\mtx{A}$.\\
\mbox{} \hfill \texttt{\color{red}B = Q'*A}
\item Form the SVD of the small matrix $\mtx{B}$: $\mtx{B} = \hat{\mtx{U}}\mtx{D}\mtx{V}^{*}$.\\
\mbox{} \hfill \texttt{\color{red}[Uhat,D,V] = svd(B,'econ')}
\item Form $\mtx{U} = \mtx{Q}\hat{\mtx{U}}$.\\
\mbox{} \hfill \texttt{\color{red}U = Q*Uhat}
\end{enumerate}
\end{minipage}}
\caption{A basic randomized algorithm. If a factorization of precisely rank $k$ is
desired, the factorization in Step 5 can be truncated to the $k$ leading terms. The
text in red is Matlab code for executing each line.}
\label{fig:RSVD}
\end{figure}

To summarize the discussion in the previous paragraph, the randomized sampling
algorithm for constructing an approximate basis for the range of a given
$m\times n$ matrix $\mtx{A}$ proceeds as follows: First pick a small integer $p$
representing how much ``over-sampling'' we do. (The choice $p=10$ is often good.)
Then execute the following steps:
\begin{enumerate}
\item Form a set of $k+p$ random Gaussian vectors $\{\vct{g}_{j}\}_{j=1}^{k+p}$.
\item Form a set $\{\vct{y}_{j}\}_{j=1}^{k+p}$ of samples from the range where
      $\vct{y}_{j} = \mtx{A}\vct{g}_{j}$.
\item Perform Gram-Schmidt on the set $\{\vct{y}_{j}\}_{j=1}^{k+p}$
      to form the ON-set $\{\vct{q}_{j}\}_{j=1}^{k+p}$.
\end{enumerate}
Now observe that the $k+p$ matrix-vector products are independent and can advantageously
be executed in parallel. A full algorithm for computing an approximate SVD using this
simplistic sampling technique for executing ``Stage A'' is summarized in Figure \ref{fig:RSVD}.

The error incurred by the randomized range finding method described in
this section is a random variable. There exist rigorous bounds for both
the expectation of this error, and for the likelihood of a large deviation
from the expectation. These bounds demonstrate that when the singular values
of $\mtx{A}$ decay ``reasonably fast,'' the error incurred is close to the
theoretically optimal one. We provide more details in Section \ref{sec:randtheory}.

\begin{remark}[How many basis vectors?]
The reader may have observed that while our stated goal was to find a matrix
$\mtx{Q}$ that holds $k$ orthonormal columns, the randomized process discussed
in this section and summarized in Figure \ref{fig:RSVD} results in a matrix
with $k+p$ columns instead. The $p$ extra vectors are needed to ensure that the
basis produced in ``Stage A'' accurately captures the $k$ dominant left singular
vectors of $\mtx{A}$. In a situation where an approximate SVD with precisely $k$
modes is sought, one can drop the last $p$ components when executing Stage B.
Using Matlab notation, we would after Step (5) run the commands
$$
\texttt{Uhat = Uhat(:,1:k); D = D(1:k,1:k); V = V(:,1:k);}.
$$
From a practical point of view, the cost of carrying around a few extra samples
in the intermediate steps is often entirely negligible.
\end{remark}

\section{Single pass algorithms}
\label{sec:singlepass}
The randomized algorithm described in Figure \ref{fig:RSVD} accesses the
matrix $\mtx{A}$ twice, first in ``Stage A'' where we build an orthonormal basis for the
column space, and then in ``Stage B'' where we project $\mtx{A}$ on to the
space spanned by the computed basis vectors. It turns out to be possible
to modify the algorithm in such a way that each entry of $\mtx{A}$ is accessed
only \textit{once.} This is important because it allows us to compute the
factorization of a matrix that is too large to be stored.

For \textit{Hermitian} matrices, the modification to Algorithm \ref{fig:RSVD}
is very minor and we describe it in Section \ref{sec:singlepass_herm}. Section
\ref{sec:singlepass_gen} then handles the case of a general matrix.

\begin{remark}[Loss of accuracy]
The single-pass algorithms described in this section tend to produce a factorization
of lower accuracy than what Algorithm \ref{fig:RSVD} would yield. In situations where
one has a choice between using either a one-pass or a two-pass algorithm, the latter is
generally preferable since it yields higher accuracy, at only moderately higher cost.
\end{remark}

\begin{remark}[Streaming Algorithms]
We say that an algorithm for processing a matrix is a \textit{streaming algorithm}
if each entry of the matrix is accessed only once, and if, in addition,
it can be fed the entries in any order. (In other words, the algorithm is not allowed
to dictate the order in which the elements are viewed.) The algorithms
described in this section satisfy both of these conditions.
\end{remark}

\subsection{Hermitian matrices}
\label{sec:singlepass_herm}

Suppose that $\mtx{A} = \mtx{A}^{*}$, and that our objective is to compute an
approximate eigenvalue decomposition
\begin{equation}
\begin{array}{cccccccccccccccccccc}
\mtx{A} &\approx& \mtx{U} & \mtx{D} & \mtx{U}^{*}\\
n\times n && n\times k & k\times k &k\times n
\end{array}
\end{equation}
with $\mtx{U}$ an orthonormal matrix and $\mtx{D}$ diagonal. (Note that for a Hermitian
matrix, the EVD and the SVD are essentially equivalent, and that the EVD is the
more natural factorization.)
Then execute Stage A with an over-sampling parameter $p$ to compute an orthonormal matrix
$\mtx{Q}$ whose columns form an approximate basis for the column space of $\mtx{A}$:
\begin{enumerate}
\item Draw a Gaussian random matrix $\mtx{G}$ of size $n\times (k+p)$.
\item Form the sampling matrix $\mtx{Y} = \mtx{A}\mtx{G}$.
\item Orthonormalize the columns of $\mtx{Y}$ to form $\mtx{Q}$, in other words $\mtx{Q} = \texttt{orth}(\mtx{Y})$.
\end{enumerate}
Then
\begin{equation}
\label{eq:dasani1}
\mtx{A} \approx \mtx{Q}\mtx{Q}^{*}\mtx{A}.
\end{equation}
Since $\mtx{A}$ is Hermitian, its row and column spaces are identical, so we also have
\begin{equation}
\label{eq:dasani2}
\mtx{A} \approx \mtx{A}\mtx{Q}\mtx{Q}^{*}.
\end{equation}
Inserting (\ref{eq:dasani1}) into (\ref{eq:dasani2}), we (informally!) find that
\begin{equation}
\label{eq:dasani3}
\mtx{A} \approx \mtx{Q}\mtx{Q}^{*}\mtx{A}\mtx{Q}\mtx{Q}^{*}.
\end{equation}
We define
\begin{equation}
\label{eq:dasani4}
\mtx{C} = \mtx{Q}^{*}\mtx{A}\mtx{Q}.
\end{equation}
If $\mtx{C}$ is known, then the post-processing is straight-forward:
Simply compute the EVD of $\mtx{C}$ to obtain $\mtx{C} = \hat{\mtx{U}}\mtx{D}\hat{\mtx{U}}^{*}$,
then define $\mtx{U} = \mtx{Q}\hat{\mtx{U}}$, to find that
\[
\mtx{A} \approx
\mtx{Q}\mtx{C}\mtx{Q}^{*} =
\mtx{Q}\hat{\mtx{U}}\mtx{D}\hat{\mtx{U}}^{*}\mtx{Q}^{*} =
\mtx{U}\mtx{D}\mtx{U}^{*}.
\]

\begin{figure}
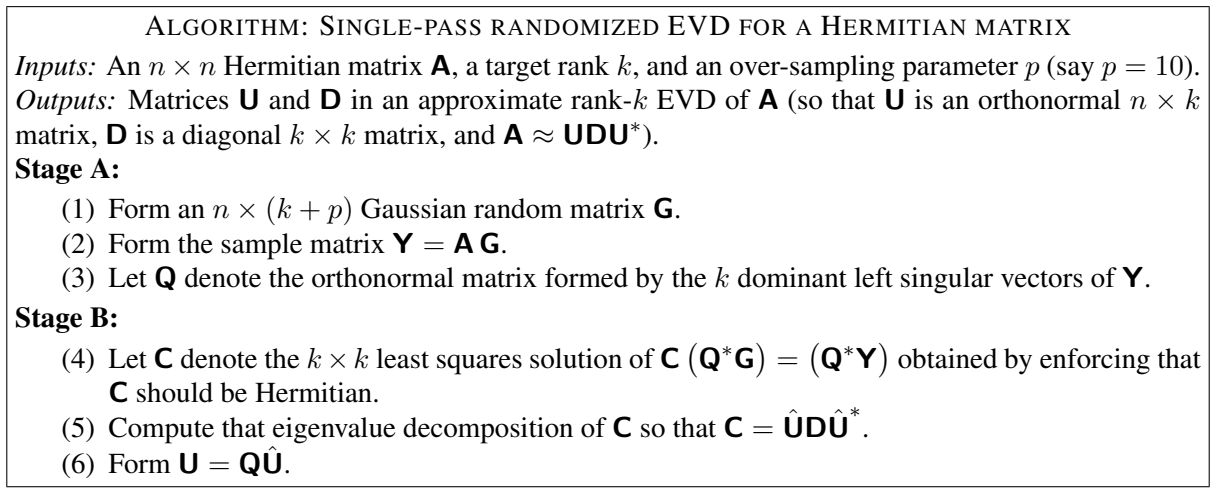

\fbox{\begin{minipage}{0.95\textwidth}
\begin{center}
\textsc{Algorithm: Single-pass randomized EVD for a Hermitian matrix}
\end{center}


\textit{Inputs:} An $n\times n$ Hermitian matrix $\mtx{A}$, a target rank $k$, and an over-sampling parameter $p$ (say $p=10$).


\textit{Outputs:} Matrices $\mtx{U}$ and $\mtx{D}$ in an approximate rank-$k$ EVD
of $\mtx{A}$ (so that $\mtx{U}$ is an orthonormal $n\times k$ matrix, $\mtx{D}$ is a diagonal $k\times k$ matrix, and
$\mtx{A} \approx \mtx{U}\mtx{D}\mtx{U}^{*}$).


\textbf{Stage A:}
\begin{enumerate}
\item Form an $n\times (k+p)$ Gaussian random matrix $\mtx{G}$.
\item Form the sample matrix $\mtx{Y} = \mtx{A}\,\mtx{G}$.
\item Let $\mtx{Q}$ denote the orthonormal matrix formed by the $k$ dominant left singular vectors of $\mtx{Y}$.
\end{enumerate}

\textbf{Stage B:}
\begin{enumerate}
\setcounter{enumi}{3}
\item Let $\mtx{C}$ denote the $k\times k$ least squares solution of
      $\mtx{C}\,\bigl(\mtx{Q}^{*}\mtx{G}\bigr)= \bigl(\mtx{Q}^{*}\mtx{Y}\bigr)$
      obtained by enforcing that $\mtx{C}$ should be Hermitian.
\item Compute that eigenvalue decomposition of $\mtx{C}$ so that
      $\mtx{C} = \hat{\mtx{U}}\mtx{D}\hat{\mtx{U}}^{*}$.
\item Form $\mtx{U} = \mtx{Q}\hat{\mtx{U}}$.
\end{enumerate}
\end{minipage}}
\caption{A basic randomized algorithm single-pass algorithm suitable for a Hermitian matrix.}
\label{fig:singlepass_herm}
\end{figure}

The problem now is that since we are seeking a single-pass algorithm, we are not in position to evaluate
$\mtx{C}$ directly from formula (\ref{eq:dasani4}). Instead, we will derive a
formula for $\mtx{C}$ that can be evaluated without revisiting $\mtx{A}$. To this
end, multiply (\ref{eq:dasani4}) by $\mtx{Q}^{*}\mtx{G}$ to obtain
\begin{equation}
\label{eq:dasani5}
\mtx{C} (\mtx{Q}^{*}\mtx{G}) = \mtx{Q}^{*}\mtx{A}\mtx{Q}\mtx{Q}^{*}\mtx{G}.
\end{equation}
We use that $\mtx{A}\mtx{Q}\mtx{Q}^{*} \approx \mtx{A}$ (cf.~(\ref{eq:dasani2})), to
approximate the right hand side in (\ref{eq:dasani5}):
\begin{equation}
\label{eq:dasani5a}
\mtx{Q}^{*}\mtx{A}\mtx{Q}\mtx{Q}^{*}\mtx{G} \approx \mtx{Q}^{*}\mtx{A}\mtx{G} = \mtx{Q}^{*}\mtx{Y}.
\end{equation}
Combining, (\ref{eq:dasani5}) and (\ref{eq:dasani5a}), and ignoring the approximation error,
we define $\mtx{C}$ as the solution of the linear system (recall that $\ell = k+p$)
\begin{equation}
\label{eq:dasani5b}
\begin{array}{ccccccccccccccccc}
\mtx{C} & \bigl(\mtx{Q}^{*}\mtx{G}\bigr) &=& \bigl(\mtx{Q}^{*}\mtx{Y}\bigr). \\
\ell \times \ell & \ell \times \ell && \ell \times \ell
\end{array}
\end{equation}
At first, it may appears that (\ref{eq:dasani5b}) is perfectly balanced in that
there are $\ell^{2}$ equations for $\ell^{2}$ unknowns. However, we need to
enforce that $\mtx{C}$ is Hermitian, so the system is actually over-determined
by roughly a factor of two. Putting everything together, we obtain the method
summarized in Figure \ref{fig:singlepass_herm}.

The procedure described in this section is less accurate than the procedure
described in Figure \ref{fig:RSVD} for two reasons: (1) The approximation error in
formula (\ref{eq:dasani3}) tends to be larger than the error in (\ref{eq:dasani1}).
\pgnotate{Follow up on this. We certainly have the bound
we have
$$
\|\mtx{A} - \mtx{P}\mtx{A}\mtx{P}\| \leq
\|\mtx{A} - \mtx{P}\mtx{A}\| + \|\mtx{P}\mtx{A} - \mtx{P}\mtx{A}\mtx{P}\| \leq
2\|\mtx{A} - \mtx{P}\mtx{A}\|.
$$
But can this worst case scenario be achieved? We also have
$$
\|\mtx{A} - \mtx{P}\mtx{A}\mtx{P}\|^{2}_{\rm Fro} \leq
\|\mtx{A} - \mtx{P}\mtx{A}\|^{2}_{\rm Fro} + \|\mtx{P}\mtx{A} - \mtx{P}\mtx{A}\mtx{P}\|^{2}_{\rm Fro}.
$$}
(2) While the matrix $\mtx{Q}^{*}\mtx{G}$ is invertible, it tends to be very ill-conditioned.

\begin{remark}[Extra over-sampling]
\label{remark:oversamp}
To combat the problem that $\mtx{Q}^{*}\mtx{G}$ tends to be ill-conditioned, it is helpful to
over-sample more aggressively when using a single pass algorithm, even to the point of setting
$p=k$ if memory allows. Once the sampling stage is completed, we form
$\mtx{Q}$ as the leading $k$ left singular vectors of $\mtx{Y}$ (compute these by forming the
full SVD of $\mtx{Y}$, and then discard the last $p$ components). Then $\mtx{C}$ will be of
size $k\times k$, and the equation that specifies $\mtx{C}$ reads
\begin{equation}
\label{eq:dasani6}
\begin{array}{ccccccccccccccccc}
\mtx{C} & \bigl(\mtx{Q}\mtx{G}\bigr) &=& \mtx{Q}^{*}\mtx{Y}.\\
k\times k & k\times \ell && k\times \ell
\end{array}
\end{equation}
Since (\ref{eq:dasani6}) is over-determined, we solve it using a least-squares technique.
Observe that we are now looking for less information (a $k\times k$ matrix rather than an $\ell\times\ell$
matrix), and have more information in order to determine it.
\end{remark}

\subsection{General matrices}
\label{sec:singlepass_gen}
\begin{figure}
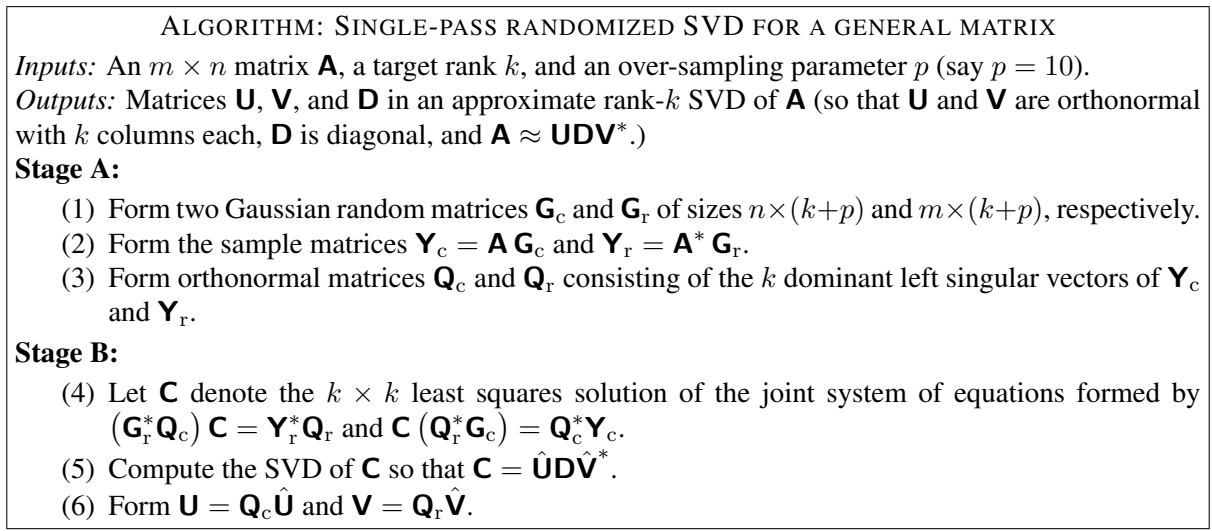

\fbox{\begin{minipage}{0.95\textwidth}
\begin{center}
\textsc{Algorithm: Single-pass randomized SVD for a general matrix}
\end{center}


\textit{Inputs:} An $m\times n$ matrix $\mtx{A}$, a target rank $k$, and an over-sampling parameter $p$ (say $p=10$).


\textit{Outputs:} Matrices $\mtx{U}$, $\mtx{V}$, and $\mtx{D}$ in an approximate rank-$k$ SVD
of $\mtx{A}$ (so that $\mtx{U}$ and $\mtx{V}$ are orthonormal with $k$ columns each, $\mtx{D}$ is diagonal, and $\mtx{A}\approx\mtx{U}\mtx{D}\mtx{V}^{*}$.)


\textbf{Stage A:}
\begin{enumerate}
\item Form two Gaussian random matrices
      $\mtx{G}_{\rm c}$ and $\mtx{G}_{\rm r}$
      of sizes $n\times (k+p)$ and $m\times (k+p)$, respectively.
\item Form the sample matrices
      $\mtx{Y}_{\rm c} = \mtx{A}    \,\mtx{G}_{\rm c}$ and
      $\mtx{Y}_{\rm r} = \mtx{A}^{*}\,\mtx{G}_{\rm r}$.
\item Form orthonormal matrices $\mtx{Q}_{\rm c}$ and $\mtx{Q}_{\rm r}$ consisting of the $k$ dominant left singular vectors
      of $\mtx{Y}_{\rm c}$ and $\mtx{Y}_{\rm r}$.
\end{enumerate}

\textbf{Stage B:}
\begin{enumerate}
\setcounter{enumi}{3}
\item Let $\mtx{C}$ denote the $k\times k$ least squares solution of the joint system
of equations formed by
$\bigl(\mtx{G}_{\rm r}^{*}\mtx{Q}_{\rm c}\bigr)\,\mtx{C} = \mtx{Y}_{\rm r}^{*}\mtx{Q}_{\rm r}$
and
$\mtx{C}\,\bigl(\mtx{Q}_{\rm r}^{*}\mtx{G}_{\rm c}\bigr) = \mtx{Q}_{\rm c}^{*}\mtx{Y}_{\rm c}$.
\item Compute the SVD of $\mtx{C}$ so that $\mtx{C} = \hat{\mtx{U}}\mtx{D}\hat{\mtx{V}}^{*}$.
\item Form $\mtx{U} = \mtx{Q}_{\rm c}\hat{\mtx{U}}$ and $\mtx{V} = \mtx{Q}_{\rm r}\hat{\mtx{V}}$.
\end{enumerate}
\end{minipage}}
\caption{A basic randomized algorithm single-pass algorithm suitable for a general matrix.}
\label{fig:singlepass_gen}
\end{figure}

We next consider a general $m\times n$ matrix $\mtx{A}$. In this case, we need to apply randomized
sampling to both its row space and its column space simultaneously. We proceed as follows:
\begin{enumerate}
\item Draw two Gaussian random matrices $\mtx{G}_{\rm c}$ of size $n\times (k+p)$ and $\mtx{G}_{\rm r}$ of size $m\times (k+p)$.
\item Form two sampling matrices $\mtx{Y}_{\rm c} = \mtx{A}\mtx{G}_{\rm c}$ and $\mtx{Y}_{\rm r} = \mtx{A}^{*}\mtx{G}_{\rm r}$.
\item Compute two basis matrices $\mtx{Q}_{\rm c} = \texttt{orth}(\mtx{Y}_{\rm c})$ and
                                 $\mtx{Q}_{\rm r} = \texttt{orth}(\mtx{Y}_{\rm r})$.
\end{enumerate}
Now define the small projected matrix via
\begin{equation}
\label{eq:dasani7}
\mtx{C} = \mtx{Q}_{\rm c}^{*}\mtx{A}\mtx{Q}_{r}.
\end{equation}
We will derive two relationships that together will determine $\mtx{C}$ in a manner that is
analogous to (\ref{eq:dasani5}). First left multiply (\ref{eq:dasani7}) by $\mtx{G}_{\rm r}^{*}\mtx{Q}_{\rm c}$ to obtain
\begin{equation}
\label{eq:dasani8}
\mtx{G}_{\rm r}^{*}\mtx{Q}_{\rm c}\mtx{C} =
\mtx{G}_{\rm r}^{*}\mtx{Q}_{\rm c}\mtx{Q}_{\rm c}^{*}\mtx{A}\mtx{Q}_{\rm r} \approx
\mtx{G}_{\rm r}^{*}\mtx{A}\mtx{Q}_{\rm r} = \mtx{Y}_{\rm r}^{*}\mtx{Q}_{\rm r}.
\end{equation}
Next we right multiply (\ref{eq:dasani7}) by $\mtx{Q}_{\rm r}^{*}\mtx{G}_{\rm c}$ to obtain
\begin{equation}
\label{eq:dasani9}
\mtx{C}\mtx{Q}_{\rm r}^{*}\mtx{G}_{\rm c} =
\mtx{Q}_{\rm c}^{*}\mtx{A}\mtx{Q}_{r}\mtx{Q}_{\rm r}^{*}\mtx{G}_{\rm c} \approx
\mtx{Q}_{\rm c}^{*}\mtx{A}\mtx{G}_{\rm c} =
\mtx{Q}_{\rm c}^{*}\mtx{Y}_{\rm c}.
\end{equation}
We now define $\mtx{C}$ as the least-square solution of the two equations
\[
\bigl(\mtx{G}_{\rm r}^{*}\mtx{Q}_{\rm c}\bigr)\,\mtx{C} = \mtx{Y}_{\rm r}^{*}\mtx{Q}_{\rm r}
\qquad\mbox{and}\qquad
\mtx{C}\,\bigl(\mtx{Q}_{\rm r}^{*}\mtx{G}_{\rm c}\bigr) = \mtx{Q}_{\rm c}^{*}\mtx{Y}_{\rm c}.
\]
Again, the system is over-determined by about a factor of 2, and it is advantageous to
make it further over-determined by more aggressive over-sampling, cf.~Remark \ref{remark:oversamp}.
Figure \ref{fig:singlepass_gen} summarizes the single-pass method for a general matrix.

\section{A method with complexity $O(mn\log k)$ for general dense matrices}
\label{sec:SRFT}
The Randomized SVD (RSVD) algorithm as given in Figure \ref{fig:RSVD}
is highly efficient when we have access to fast methods for evaluating matrix-vector
products $\vct{x} \mapsto \mtx{A}\vct{x}$. For the case where $\mtx{A}$
is a general $m\times n$ matrix given simply as an array or real numbers,
the cost of evaluating the sample matrix $\mtx{Y} = \mtx{A}\mtx{G}$
(in Step (2) of the algorithm in Figure \ref{fig:RSVD}) is $O(mnk)$.
RSVD is still often faster than classical methods
since the matrix-matrix multiply can be highly optimized, but it does
not have an edge in terms of asymptotic complexity. However, it turns out
to be possible to modify the algorithm by replacing the Gaussian random
matrix $\mtx{G}$ with a different random matrix $\mtx{\Omega}$ that has
two seemingly contradictory properties:

\vspace{1mm}

\noindent
(1) $\mtx{\Omega}$ is sufficiently \textit{structured} that
$\mtx{A}\mtx{\Omega}$ can be evaluated in $O(mn\log(k))$ flops;

\vspace{1mm}

\noindent
(2) $\mtx{\Omega}$\! is\! sufficiently\! \textit{random} that the columns of
$\mtx{A}\mtx{\Omega}$ accurately span the range of $\mtx{A}$.

\vspace{1mm}

\noindent
For instance, a good choice of random matrix $\mtx{\Omega}$ is
\begin{equation}
\label{eq8:SRFT}
\begin{array}{ccccccccccccccc}
\mtx{\Omega} &=& \mtx{D}&\mtx{F}&\mtx{S},\\
n\times \ell && n\times n & n \times n & n\times \ell
\end{array}
\end{equation}
where $\mtx{D}$ is a diagonal matrix whose diagonal entries are complex
numbers of modulus one drawn from a uniform distribution on the unit circle
in the complex plane, where $\mtx{F}$ is the discrete Fourier transform,
\[
\mtx{F}(p,q) = n^{-1/2}\,e^{-2\pi\mathrm{i}(p-1)(q-1)/n},
\qquad p,q \in \{1,\,2,\,3,\,\dots,\,n\},
\]
and where $\mtx{S}$ is a matrix consisting of a random subset of $\ell$
columns from the $n\times n$ unit matrix (drawn without replacement).
In other words, given an arbitrary matrix $\mtx{X}$ of size $m\times n$,
the matrix $\mtx{X}\mtx{S}$ consists of a randomly drawn subset of $\ell$
columns of $\mtx{X}$.
For the matrix $\mtx{\Omega}$ specified by (\ref{eq8:SRFT}), the product
$\mtx{X}\mtx{\Omega}$ can be evaluated via a subsampled FFT in $O(mn\log(\ell))$
operations. The parameter $\ell$ should be chosen slightly larger than the
target rank $k$; the choice $\ell = 2k$ is often good. (A transform of this
type was introduced in \cite{2006_ailon_chazelle_FJLT} under the name
``Fast Johnson-Lindenstrauss Transform'' and was applied to the problem of
low-rank approximation in \cite{2007_woolfe_liberty_rokhlin_tygert,2006_sarlos_improved}. See
also \cite{2013_ailon_almost,2014_kane_nelson_sparser_JL,2009_edo_liberty_dissertation}.)

By using the structured random matrix described in this section, we can
reduce the complexity of ``Stage A'' in the RSVD from $O(mnk)$ to $O(mn\log (k))$.
In order to attain overall cost $O(mn\log (k))$, we must also modify ``Stage B''
to eliminate the need to compute $\mtx{Q}^{*}\mtx{A}$ (since direct evaluation of
$\mtx{Q}^{*}\mtx{A}$ has cost $O(mnk)$).
One option is to use the single pass algorithm described in \ref{fig:singlepass_gen},
using the structured random matrix to approximate both the row and the column spaces
of $\mtx{A}$. A second, and typically better, option is to use a so called
\textit{row-extraction} technique for Stage B; we describe the details in
Section \ref{sec:ID}.

Our theoretical understanding of the errors incurred by the accelerated range finder
is not as satisfactory as what we have for Gaussian random matrices,
cf.~\cite{2011_martinsson_randomsurvey}*{Sec.~11}. In the general case,
only quite weak results have been proven. In practice, the accelerated scheme is often as
accurate as the Gaussian one, but we do not currently have good theory to predict precisely
when this happens.

\section{Theoretical performance bounds}
\label{sec:randtheory}

In this section, we will briefly summarize some proven results concerning the
error in the output of the basic RSVD algorithm in Figure \ref{fig:RSVD}. Observe
that the factors $\mtx{U},\,\mtx{D},\,\mtx{V}$ depend not only on $\mtx{A}$, but
also on the draw of the random matrix $\mtx{G}$. This means that the error that
we try to bound is a \textit{random variable.} It is therefore natural to seek
bounds on first the expected value of the error, and then on the
likelihood of large deviations from the expectation.

Before we start, let us recall from Remark \ref{remark8:nomoreerror} that all the
error incurred by the RSVD algorithm in Figure \ref{fig:RSVD} is incurred in Stage A.
The reason is that the ``post-processing'' in Stage B is  exact (up to floating point
arithmetic). Consequently, we can (and will) restrict ourselves to providing bounds
on $\|\mtx{A} - \mtx{Q}\mtx{Q}^{*}\mtx{A}\|$.

\begin{remark}
The theoretical investigation of errors resulting from randomized methods in linear
algebra is an active area of research that draws heavily on random matrix theory,
theoretical computer science, classical numerical linear algebra, and many other fields.
Our objective here is merely to state a couple of representative results, without providing
any proofs or details about their derivation. Both results are taken from
\cite{2011_martinsson_randomsurvey}, where the interested reader can find
an in-depth treatment of the subject. More recent results pertaining to
the RSVD can be found in, e.g., \cite{2015_candes_rsvd_bounds,2015_gu_randomized_subspaceiteration}.
\end{remark}

\subsection{Bounds on the expectation of the error}
A basic result on the \textit{typical} error observed is
Theorem 10.6 of \cite{2011_martinsson_randomsurvey}, which states:

\begin{thm}
\label{thm:RSVDexpectationerror}
Let $\mtx{A}$ be an $m\times n$ matrix with singular values $\{\sigma_{j}\}_{j=1}^{\min(m,n)}$.
Let $k$ be a target rank, and let $p$ be an over-sampling parameter
such that $p \geq 2$ and $k+p \leq \min(m,n)$. Let $\mtx{G}$ be a
Gaussian random matrix of size $n\times (k+p)$ and set
{\rm $\mtx{Q} = \texttt{orth}(\mtx{A}\mtx{G})$}.
Then the average error, as measured in the Frobenius norm, satisfies
\begin{equation}
\label{eq8:expectationbound}
\mathbb{E}\bigl[\|\mtx{A} - \mtx{Q}\mtx{Q}^{*}\mtx{A}\|_{\rm Fro}\bigr] \leq
\left(1 + \frac{k}{p-1}\right)^{1/2}\,
\left(\sum_{j=k+1}^{\min(m,n)} \sigma_{j}^{2}\right)^{1/2},
\end{equation}
where $\mathbb{E}$ refers to expectation with respect to the draw of $\mtx{G}$.
The corresponding result for the spectral norm reads
\begin{equation}
\label{eq:expectationboundspectral}
\mathbb{E}\bigl[ \|\mtx{A} - \mtx{Q}\mtx{Q}^{*}\mtx{A}\|\bigr]
    \leq \left(1 + \sqrt{\frac{k}{p-1}} \right) \sigma_{k+1}
        + \frac{e\sqrt{k+p}}{p}
        \left(\sum_{j=k+1}^{\min(m,n)} \sigma_{j}^2 \right)^{1/2}.
\end{equation}
\end{thm}

When errors are measured in the \textit{Frobenius norm,} Theorem \ref{thm:RSVDexpectationerror}
is very gratifying. For our standard recommendation of $p=10$, we are basically within a factor
of $\sqrt{1+k/9}$ of the theoretically minimal error. (Recall that the Eckart-Young theorem
states that $\left(\sum_{j=k+1}^{\min(m,n)}\sigma_{j}^{2}\right)^{1/2}$ is a lower bound on the
residual for any rank-$k$ approximant.) If you over-sample more aggressively and set
$p=k+1$, then we are within a distance of $\sqrt{2}$ of the theoretically minimal error.

When errors are measured in the \textit{spectral norm,} the situation is much less rosy. The
first term in the bound in (\ref{eq:expectationboundspectral}) is perfectly acceptable, but
the second term is unfortunate in that it involves the minimal error in the Frobenius norm,
which can be much larger, especially when $m$ or $n$ are large. The
theorem is quite sharp, as it turns out, so the sub-optimality expressed in (\ref{eq:expectationboundspectral})
reflects a true limitation on the accuracy to be expected from the basic randomized scheme.

The extent to which the error in (\ref{eq:expectationboundspectral}) is problematic
depends on how rapidly the ``tail'' singular values $\{\sigma_{j}\}_{j=1}^{\min(m,n)}$ decay. If they
decay fast, then the spectral norm error and the Frobenius norm error are similar, and the
RSVD works well. If they decay slowly, then the RSVD performs fine when errors are measured in
the Frobenius norm, but not very well when the spectral norm is the one of interest. To illustrate
the difference, let us consider two situations:

\begin{widequote}
\textit{Case 1 --- fast decay:} Suppose that the tail singular values decay exponentially fast,
so that for some $\beta \in (0,1)$ we have $\sigma_{j} \approx \sigma_{k+1}\,\beta^{j-k-1}$ for $j > k$. Then
$\left(\sum_{j=k+1}^{\min(m,n)}\sigma_{j}^{2}\right)^{1/2} \approx
\sigma_{k+1}\left(\sum_{j=k+1}^{\min(m,n)}\beta^{2(j-k-1)}\right)^{1/2} \leq
\sigma_{k+1}(1 - \beta^{2})^{-1/2}$. As long as $\beta$ is not very close to 1, we see that the
contribution from the tail singular values is modest in this case.


\noindent
\textit{Case 2 --- no decay:} Suppose that the tail singular values exhibit \textit{no} decay,
so that $\sigma_{j} = \sigma_{k+1}$ for $j>k$. This represents the worst case scenario, and now
$\left(\sum_{j=k+1}^{\min(m,n)}\sigma_{j}^{2}\right)^{1/2} = \sigma_{k+1}\sqrt{\min(m,n)-k}$. Since we
want to allow for $n$ and $m$ to be very large, this represents devastating suboptimality.
\end{widequote}

Fortunately, it is possible to modify the RSVD in such a way that the errors produced
are close to optimal in both the spectral and the Frobenius norms. The price to pay is
a modest increase in the computational cost. See Section \ref{sec:power} and
\cite{2011_martinsson_randomsurvey}*{Sec.~4.5}.


\subsection{Bounds on the likelihood of large deviations}
One can prove that (perhaps surprisingly) the likelihood of a large deviation
from the mean depends only on the over-sampling parameter $p$, and decays
extraordinarily fast. For instance, one can prove that if $p \geq 4$, then
\begin{equation}
\label{eq8:tailbound}
||\mtx{A} - \mtx{Q}\mtx{Q}^{*}\mtx{A}||
\leq \left( 1 + 17 \sqrt{1 + k/p} \right) \sigma_{k+1}
+ \frac{8\sqrt{k+p}}{p+1} \left( \sum\nolimits_{j > k} \sigma_j^2 \right)^{1/2},
\end{equation}
with failure probability at most $3\,e^{-p}$, see \cite{2011_martinsson_randomsurvey}*{Cor.~10.9}.

\section{An accuracy enhanced randomized scheme}
\label{sec:power}

\subsection{The key idea --- power iteration}
\label{sec:poweridea}
We saw in Section \ref{sec:randtheory} that the basic randomized scheme (see, e.g.,
Figure \ref{fig:RSVD}) gives accurate results for matrices whose
singular values decay rapidly, but tends to produce suboptimal
results when they do not. To recap,
suppose that we compute a rank-$k$ approximation to an $m\times n$
matrix $\mtx{A}$ with singular values $\{\sigma_{j}\}_{j=1}^{\min(m,n)}$. The theory
shows that the error measured in the spectral norm is bounded only by
a factor that scales with
$\bigl(\sum_{j>k}\sigma_{j}^{2}\bigr)^{1/2}$. When the singular
values decay slowly, this quantity can be much larger than the
theoretically minimal approximation error (which is $\sigma_{k+1}$).

Recall that the objective of the randomized sampling is to construct
a set of orthonormal vectors $\{\vct{q}_{j}\}_{j=1}^{\ell}$ that capture to
high accuracy the space spanned by the $k$ dominant left singular
vectors $\{\vct{u}_{j}\}_{j=1}^{k}$ of $\mtx{A}$. The idea is now
to sample not $\mtx{A}$, but the matrix $\mtx{A}^{(q)}$ defined by
\[
\mtx{A}^{(q)} = \bigl(\mtx{A}\mtx{A}^{*}\bigr)^{q}\mtx{A},
\]
where $q$ is a small positive integer (say, $q=1$ or $q=2$).
A simple calculation shows that if $\mtx{A}$ has the SVD
$\mtx{A} = \mtx{U}\mtx{D}\mtx{V}^{*}$, then the SVD of $\mtx{A}^{(q)}$ is
\[
\mtx{A}^{(q)} = \mtx{U}\,\mtx{D}^{2q+1}\,\mtx{V}^{*}.
\]
In other words, $\mtx{A}^{(q)}$ has the same left singular vectors
as $\mtx{A}$, while its singular values are $\{\sigma_{j}^{2q+1}\}_{j}$.
Even when the singular values of $\mtx{A}$ decay slowly, the
singular values of $\mtx{A}^{(q)}$ tend to decay fast enough for
our purposes.

\begin{figure}
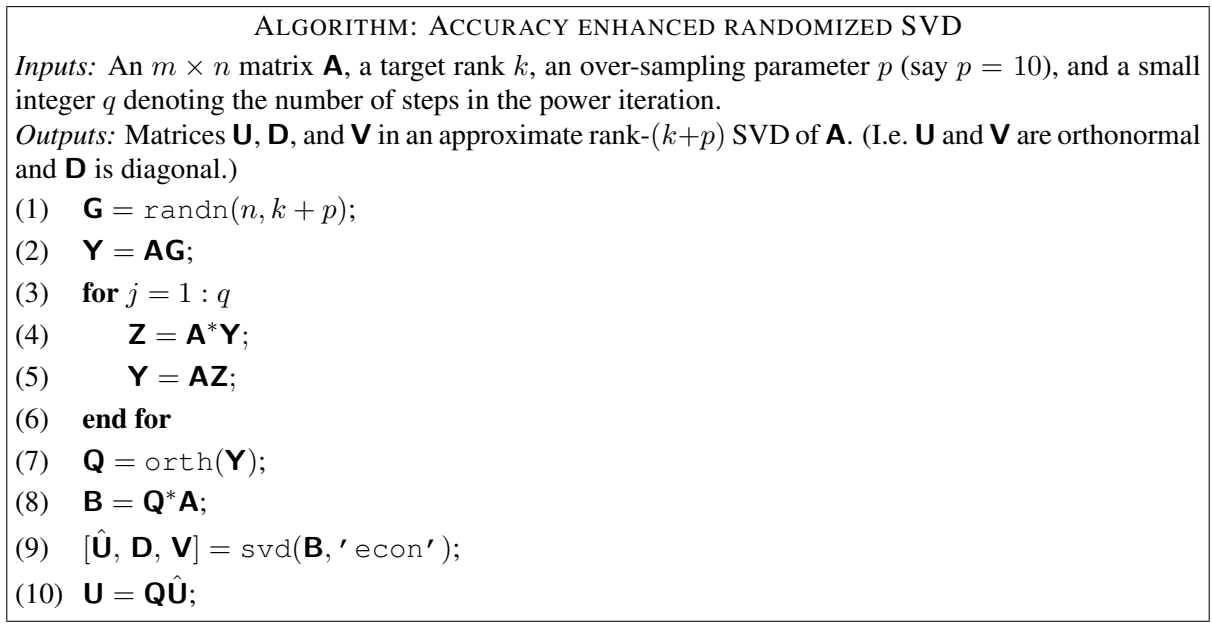

\fbox{\begin{minipage}{0.95\textwidth}
\begin{center}
\textsc{Algorithm: Accuracy enhanced randomized SVD}
\end{center}


\textit{Inputs:} An $m\times n$ matrix $\mtx{A}$, a target rank $k$, an over-sampling parameter $p$ (say $p=10$),
and a small integer $q$ denoting the number of steps in the power iteration.


\textit{Outputs:} Matrices $\mtx{U}$, $\mtx{D}$, and $\mtx{V}$ in an approximate rank-$(k+p)$ SVD
of $\mtx{A}$. (I.e.~$\mtx{U}$ and $\mtx{V}$ are orthonormal and $\mtx{D}$ is diagonal.)


\begin{tabbing}
\hspace{8mm} \= \hspace{5mm} \= \hspace{5mm} \= \hspace{15mm}\kill
(1) \> $\mtx{G} = \texttt{randn}(n,k+p)$;\\[1mm]
(2) \> $\mtx{Y} = \mtx{A}\mtx{G}$;\\[1mm]
(3) \> \textbf{for} $j = 1:q$\\[1mm]
(4) \> \> $\mtx{Z} = \mtx{A}^{*}\mtx{Y};$\\[1mm]
(5) \> \> $\mtx{Y} = \mtx{A}\mtx{Z};$\\[1mm]
(6) \> \textbf{end for}\\[1mm]
(7) \> $\mtx{Q} = \texttt{orth}(\mtx{Y})$;\\[1mm]
(8) \> $\mtx{B} = \mtx{Q}^{*}\mtx{A}$;\\[1mm]
(9) \> $[\hat{\mtx{U}},\,\mtx{D},\,\mtx{V}] = \texttt{svd}(\mtx{B},\texttt{'econ'})$;\\[1mm]
(10) \> $\mtx{U} = \mtx{Q}\hat{\mtx{U}}$;
\end{tabbing}

\end{minipage}}
\caption{The accuracy enhanced randomized SVD. If a factorization of precisely rank $k$ is
desired, the factorization in Step 9 can be truncated to the $k$ leading terms.}
\label{fig:RSVDP}
\end{figure}

The accuracy enhanced scheme now consists of drawing a Gaussian matrix $\mtx{G}$ and
then forming a sample matrix
\[
\mtx{Y} = \bigl(\mtx{A}\mtx{A}^{*}\bigr)^{q}\mtx{A}\mtx{G}.
\]
Then orthonormalize the columns of $\mtx{Y}$ to obtain
$\mtx{Q} = \texttt{orth}(\mtx{Y})$, and proceed as before.
The resulting scheme is shown in Figure \ref{fig:RSVDP}.

\begin{remark}
\label{remark:RSVDP_stabilized}
The scheme described in Figure \ref{fig:RSVDP} can lose accuracy due to round-off
errors. The problem is that as $q$ increases, all columns in the sample matrix
$\mtx{Y} = \bigl(\mtx{A}\mtx{A}^{*}\bigr)^{q}\mtx{A}\mtx{G}$ tend to align closer
and closer to the dominant left singular vector. This means that essentially all
information about the singular values and singular vectors associated with smaller
singular values get lots to round-off errors. Roughly speaking, if
\[
\frac{\sigma_{j}}{\sigma_{1}} \leq \epsilon_{\rm mach}^{1/(2q+1)},
\]
where $\epsilon_{\rm mach}$ is machine precision,
then all information associated with the $j$'th singular value and beyond is lost
(see Section 3.2 of \cite{2016_martinsson_randomized_notes}).
This problem can be fixed by orthonormalizing the columns
between each iteration, as shown in Figure \ref{fig:RSVDP_stabilized}.
The modified scheme is more costly due to the extra calls to \texttt{orth}.
(However, note that \texttt{orth} can be executed using
\textit{unpivoted} Gram-Schmidt, which is quite fast.)
\end{remark}

\begin{figure}
\fbox{\begin{minipage}{0.95\textwidth}
\begin{center}
\textsc{Algorithm: Accuracy enhanced randomized SVD\\ (with orthonormalization)}
\end{center}

\begin{tabbing}
\hspace{8mm} \= \hspace{5mm} \= \hspace{5mm} \= \hspace{15mm}\kill
(1) \> $\mtx{G} = \texttt{randn}(n,k+p)$;\\[1mm]
(2) \> $\mtx{Q} = \texttt{orth}(\mtx{A}\mtx{G})$;\\[1mm]
(3) \> \textbf{for} $j = 1:q$\\[1mm]
(4) \> \> $\mtx{W} = \texttt{orth}(\mtx{A}^{*}\mtx{Q});$\\[1mm]
(5) \> \> $\mtx{Q} = \texttt{orth}(\mtx{A}\mtx{W});$\\[1mm]
(6) \> \textbf{end for}\\[1mm]
(7) \> $\mtx{B} = \mtx{Q}^{*}\mtx{A}$;\\[1mm]
(8) \> $[\hat{\mtx{U}},\,\mtx{D},\,\mtx{V}] = \texttt{svd}(\mtx{B},\texttt{'econ'})$;\\[1mm]
(9) \> $\mtx{U} = \mtx{Q}\hat{\mtx{U}}$;
\end{tabbing}

\end{minipage}}
\caption{This algorithm takes the same inputs and outputs as the method in Figure \ref{fig:RSVDP}.
The only difference is that orthonormalization is carried out between each step of the power iteration,
to avoid loss of accuracy due to rounding errors.}
\label{fig:RSVDP_stabilized}
\end{figure}

\subsection{Theoretical results}
A detailed error analysis of the scheme described in Figure \ref{fig:RSVDP} is provided
in \cite{2011_martinsson_randomsurvey}*{Sec.~10.4}. In particular, the key theorem states:

\begin{thm}
Let $\mtx{A}$ denote an $m\times n$ matrix, let $p \geq 2$ be an over-sampling parameter, and let
$q$ denote a small integer. Draw a Gaussian matrix $\mtx{G}$ of size $n\times (k+p)$, set
$\mtx{Y} = (\mtx{A}\mtx{A}^{*})^{q} \mtx{A}\mtx{G}$, and let $\mtx{Q}$ denote an $m\times (k+p)$
orthonormal matrix
resulting from orthonormalizing the columns of $\mtx{Y}$.
Then
\begin{equation}
\label{eq:averageerrorpower}
\mathbb{E}\bigl[\|\mtx{A} - \mtx{Q}\mtx{Q}^{*}\mtx{A}\|\bigr]
    \leq \left[ \left( 1 + \sqrt{\frac{k}{p-1}} \right) \sigma_{k+1}^{2q+1}
    + \frac{e\sqrt{k+p}}{p} \left( \sum_{j=k+1}^{\min(m,n)} \sigma_j^{2(2q+1)} \right)^{1/2} \right]^{1/(2q+1)}.
\end{equation}
\end{thm}

The bound in (\ref{eq:averageerrorpower}) is slightly opaque. To simplify it, let us
consider a worst case scenario where there is no decay in the singular values beyond the
truncation point, so that $\sigma_{k+1} = \sigma_{k+2} = \cdots = \sigma_{\min(m,n)}$. Then
(\ref{eq:averageerrorpower}) simplifies to
\[
\mathbb{E}\bigl[\|\mtx{A} - \mtx{Q}\mtx{Q}^{*}\mtx{A}\|\bigr]
    \leq \left[ 1 + \sqrt{\frac{k}{p-1}}
    + \frac{e\sqrt{k+p}}{p} \cdot \sqrt{ \min\{m,n\} - k } \right]^{1/(2q+1)}
    \sigma_{k+1}.
\]
In other words, as we increase the exponent $q$, the power scheme drives the factor
that multiplies $\sigma_{k+1}$ to one exponentially fast. This factor represents the
degree of ``sub-optimality'' you can expect to see.

\subsection{Extended sampling matrix}
\label{sec:extendedsamplematrix}
The scheme described in Section \ref{sec:poweridea} is slightly wasteful in that
it does not directly use all the sampling vectors computed. To further improve accuracy,
one can for a small positive integer $q$ form an ``extended'' sampling matrix
\[
\mtx{Y} = \bigl[\mtx{A}\mtx{G},\,\mtx{A}^{2}\mtx{G},\,\dots,\,\mtx{A}^{q}\mtx{G}\bigr].
\]
Observe that this new sampling matrix $\mtx{Y}$ has $q\ell$ columns. Then proceed as before:
\begin{equation}
\label{eq:red}
\mtx{Q} = \texttt{qr}(\mtx{Y}),\qquad
\mtx{B} = \mtx{Q}^{*}\mtx{A},\qquad
[\hat{\mtx{U}},\mtx{D},\mtx{V}] = \texttt{svd}(\mtx{B},\texttt{'econ'}),\qquad
\mtx{U} = \mtx{Q}\hat{\mtx{U}}.
\end{equation}
The computations in (\ref{eq:red}) can be quite expensive since the ``tall thin''
matrices being operated on now have $q\ell$ columns, rather than the tall thin
matrices in, e.g., the algorithm in Figure \ref{fig:RSVDP}, which have only $\ell$
columns. This results in an increase in cost for all operations (QR factorization,
matrix-matrix multiply, SVD) by a factor of $O(q^{2})$. Consequently, the scheme
described here is primarily useful in situations where the computational cost is
dominated by applications of $\mtx{A}$ and $\mtx{A}^{*}$, and we want to maximally
leverage all interactions with $\mtx{A}$. An early discussion of this idea
can be found in \cite[Sec.~4.4]{2009_szlam_power}, with a more detailed discussion
in \cite{2015_musco_randomized_block_krylov}.

\section{The Nystr\"om method for symmetric positive definite matrices}

When the input matrix $\mtx{A}$ is symmetric positive definite (spd), the
\textit{Nystr\"om method} can be used to improve
the quality of standard factorizations at almost no additional cost;
see~\cite{2005_drineas_nystrom} and its bibliography.
To describe the idea, we first recall from Section \ref{sec:singlepass_herm}
that when $\mtx{A}$ is Hermitian (which of course every spd matrix is),
then it is natural to use the approximation
\begin{equation}
\label{eq:spd_standard}
\mtx{A} \approx \mtx{Q}\,\bigl(\mtx{Q}^{*}\mtx{A}\mtx{Q}\bigr)\,\mtx{Q}^{*}.
\end{equation}
In contrast, the so called ``Nystr\"{o}m scheme'' relies on the rank-$k$ approximation
\begin{equation}
\label{eq:spd_nystrom}
\mtx{A}
\approx (\mtx{A}\mtx{Q})\,\bigl(\mtx{Q}^{*}\mtx{A}\mtx{Q}\bigr)^{-1}\,(\mtx{A}\mtx{Q})^{*}.
\end{equation}
For both stability and computational efficiency, we typically rewrite (\ref{eq:spd_nystrom}) as
\[
\mtx{A} \approx \mtx{F}\mtx{F}^{*},
\]
where $\mtx{F}$ is an approximate \textit{Cholesky} factor of $\mtx{A}$ of size $n\times k$, defined by
\[
\mtx{F} = (\mtx{A}\mtx{Q})\,\bigl(\mtx{Q}^{*}\mtx{A}\mtx{Q}\bigr)^{-1/2}.
\]
To compute the factor $\mtx{F}$ numerically, first form the matrices
$\mtx{B}_1 = \mtx{AQ}$ and $\mtx{B}_2 = \mtx{Q}^{*} \mtx{B}_1$.
Observe that $\mtx{B}_{2}$ is necessarily spd, which means that we
can compute its Cholesky factorization $\mtx{B}_2 = \mtx{C}^{*} \mtx{C}$.
Finally compute the factor $\mtx{F} = \mtx{B}_1 \mtx{C}^{-1}$
by performing a triangular solve. The low-rank factorization~\eqref{eq:spd_nystrom} can be
converted to a standard decomposition using the techniques from Section \ref{sec:twostage}.

The Nystr{\"o}m technique for computing an approximate eigenvalue decomposition is
given in Figure \ref{fig:nystrom}. Let us compare the cost of this method to the more
straight-forward method resulting from using the formula (\ref{eq:spd_standard}). In
both cases, we need to twice apply $\mtx{A}$ to a set of $k+p$ vectors (first in
computing $\mtx{A}\mtx{G}$, then in computing $\mtx{A}\mtx{Q})$. But the Nystr\"om
method tends to result in substantially more accurate results. Informally speaking,
the reason is that by exploiting the spd property of $\mtx{A}$, we can take one step
of power iteration ``for free.''
For a more formal analysis of the cost and accuracy of the Nystr\"om method, we refer
the reader to \cite{2016_gittens_revisiting_nystrom,2016_becker_nystrom}.

\begin{figure}[tb]
\begin{center}
\fbox{
\begin{minipage}{.9\textwidth}
\begin{center}
\textsc{Algorithm: Eigenvalue Decomposition via the Nystr{\"o}m Method}
\end{center}

\textit{Given an $n\times n$ non-negative matrix $\mtx{A}$, a target rank $k$ and an over-sampling parameter $p$,
this procedure computes an approximate eigenvalue decomposition $\mtx{A} \approx \mtx{U}\mtx{\Lambda}\mtx{U}^{*}$,
where $\mtx{U}$ is orthonormal, and $\mtx{\Lambda}$ is nonnegative and diagonal.}

\begin{enumerate}
\item Draw a Gaussian random matrix $\mtx{G} = \texttt{randn}(n,k+p)$.
\item Form the sample matrix $\mtx{Y} = \mtx{A}\mtx{G}$.
\item Orthonormalize the columns of the sample matrix to obtain the basis matrix $\mtx{Q} = \texttt{orth}(\mtx{Y})$.
\item Form the matrices $\mtx{B}_{1} = \mtx{A}\mtx{Q}$ and $\mtx{B}_{2} = \mtx{Q}^{*}\mtx{B}_{1}$.
\item Perform a Cholesky factorization $\mtx{B}_{2} = \mtx{C}^{*} \mtx{C}$.
\item Form $\mtx{F} = \mtx{B}_{1}\mtx{C}^{-1}$ using a triangular solve.
\item Compute an SVD of the Cholesky factor $[\mtx{U},\,\mtx{\Sigma},\sim] = \texttt{svd}(\mtx{F},\texttt{'econ'})$.
\item Set $\mtx{\Lambda} = \mtx{\Sigma}^{2}$.
\end{enumerate}
\end{minipage}}
\end{center}
\caption{The Nystr\"om method for low-rank approximation of self-adjoint matrices
with non-negative eigenvalues. It involves two applications of $\mtx{A}$ to matrices with $k+p$
columns, and has comparable cost to the basic RSVD in Figure \ref{fig:RSVD}.
However, it exploits the symmetry of $\mtx{A}$ to boost the accuracy.}
\label{fig:nystrom}
\end{figure}

\section{Randomized algorithms for computing the Interpolatory Decomposition (ID)}
\label{sec:ID}

\subsection{Structure preserving factorizations}
\label{sec:IDbasic}
Any matrix $\mtx{A}$ of size $m\times n$ and rank $k$, where $k < \min(m,n)$, admits a
so called ``interpolative decomposition (ID)'' which takes the form
\begin{equation}
\label{eq:defID1}
\begin{array}{cccc}
\mtx{A} &=& \mtx{C} &\mtx{Z},\\
m\times n && m\times k & k\times n
\end{array}
\end{equation}
where the matrix $\mtx{C}$ is given by a subset of the columns of $\mtx{A}$
and where $\mtx{Z}$ is well-conditioned in a sense that we will make precise shortly.
The ID has several advantages, as compared to, e.g., the QR or SVD factorizations:
\begin{itemize}
\item If $\mtx{A}$ is sparse or non-negative, then $\mtx{C}$ shares these properties.
\item The ID requires less memory to store than either the QR or the singular value decomposition.
\item Finding the indices associated with the spanning columns is often helpful in
\textit{data interpretation.}
\item In the context of numerical algorithms for discretizing PDEs and integral equations,
      the ID often preserves ``the physics'' of a problem in a way that the QR or SVD do not.
\end{itemize}
One shortcoming of the ID is that when $\mtx{A}$ is not of precisely rank $k$, then the
approximation error by the best possible rank-$k$ ID can be substantially larger than
the theoretically minimal error. (In fact, the ID and the column pivoted QR factorizations
are closely related, and they attain \textit{exactly} the same minimal error.)

For future reference, let $J_{\rm s}$ be an index vector in $\{1,2,\dots,n\}$ that
identifies the $k$ columns in $\mtx{C}$ so that
$$
\mtx{C} = \mtx{A}(:,J_{\rm s}).
$$
One can easily show (see, e.g., \cite{2016_martinsson_randomized_notes}*{Thm.~9}) that
any matrix of rank $k$ admits a factorization (\ref{eq:defID1}) that is well-conditioned
in the sense that each entry of $\mtx{Z}$ is bounded in modulus by one. However, any
algorithm that is guaranteed to find such an optimally conditioned factorization must
have combinatorial complexity. Polynomial time algorithms with high practical efficiency
are discussed in \cite{gu1996,2005_martinsson_skel}. Randomized algorithms are described
in \cite{2011_martinsson_randomsurvey,2014_martinsson_CUR}.

\begin{remark}
The interpolative decomposition is closely related to the so called CUR decomposition
which has been studied extensively in the context of randomized algorithms
\cite{2009_mahoney_CUR,2013_shusen_CUR,boutsidis2009improved,drineas2012fast}. We will
return to this point in Section \ref{sec:CUR}.
\end{remark}

\subsection{Three flavors of ID: The row, column, and double-sided ID}

Section \ref{sec:IDbasic} describes a factorization
where we use a subset of the \textit{columns} of $\mtx{A}$ to span its \textit{column space}.
Naturally, this factorization has a sibling which uses the \textit{rows} of $\mtx{A}$ to span its
\textit{row space.} In other words $\mtx{A}$ also admits the factorization
\begin{equation}
\label{eq:defID2}
\begin{array}{cccc}
\mtx{A} &=& \mtx{X} &\mtx{R},\\
m\times n && m\times k & k\times n
\end{array}
\end{equation}
where $\mtx{R}$ is a matrix consisting of $k$ rows of $\mtx{A}$, and where
$\mtx{X}$ is a matrix that contains the $k\times k$ identity matrix. We let
$I_{\rm s}$ denote the index vector of length $k$ that marks the ``skeleton'' rows
so that $\mtx{R} = \mtx{A}(I_{\rm s},:)$.

Finally, there exists a so called \textit{double-sided ID} which takes the form
\begin{equation}
\label{eq:defID3}
\begin{array}{ccccc}
\mtx{A} &=& \mtx{X} & \mtx{A}_{\rm s} & \mtx{Z},\\
m\times n && m\times k & k\times k & k\times n
\end{array}
\end{equation}
where $\mtx{X}$ and $\mtx{Z}$ are the same matrices as those that appear in
(\ref{eq:defID1}) and (\ref{eq:defID2}), and where $\mtx{A}_{\rm s}$ is the
$k\times k$ submatrix of $\mtx{A}$ given by
$$
\mtx{A}_{\rm s} = \mtx{A}(I_{\rm s},J_{\rm s}).
$$

\subsection{Deterministic techniques for computing the ID}
In this section we demonstrate that there is a close connection
between the column ID and the classical column pivoted QR factorization
(CPQR). The end result is that standard software used to compute the CPQR
can with some light post-processing be used to compute the column ID.

As a starting point, recall that for a given $m\times n$ matrix $\mtx{A}$,
with $m \geq n$, the QR factorization can be written as
\begin{equation}
\label{eq:defCPQS}
\begin{array}{ccccccccccccccccccc}
\mtx{A}&\mtx{P} &=& \mtx{Q}&\mtx{S},\\
m\times n & n\times n && m\times n & n\times n
\end{array}
\end{equation}
where $\mtx{P}$ is a permutation matrix, where $\mtx{Q}$ has orthonormal
columns and where $\mtx{S}$ is upper triangular. \footnote{We use the letter $\mtx{S}$
instead of the traditional $\mtx{R}$ to avoid confusion with the ``R''-factor
in the row ID, (\ref{eq:defID2}).} Since our objective here is to construct
a rank-$k$ approximation to $\mtx{A}$, we split off the leading $k$ columns
from $\mtx{Q}$ and $\mtx{S}$ to obtain partitions
\begin{equation}
\label{eq:CPQSpart}
\mtx{Q} =
\kbordermatrix{&
k & n-k\\
m &\mtx{Q}_{1} & \mtx{Q}_{2}},
\qquad\mbox{and}\qquad
\mtx{S} =
\kbordermatrix{&
k & n-k\\
k &\mtx{S}_{11} & \mtx{S}_{12}\\
m-k &\mtx{0} &\mtx{S}_{22}}.
\end{equation}
Combining (\ref{eq:defCPQS}) and (\ref{eq:CPQSpart}), we then find that
\begin{equation}
\label{eq:CPQSidentity}
\mtx{A}\mtx{P} =
[\mtx{Q}_{1}\ \ |\ \ \mtx{Q}_{2}]
\left[\begin{array}{cc}
\mtx{S}_{11} & \mtx{S}_{12} \\
\mtx{0}      & \mtx{S}_{22}
\end{array}\right] =
[\mtx{Q}_{1}\mtx{S}_{11}\ \ |\ \ \mtx{Q}_{1}\mtx{S}_{12} + \mtx{Q}_{2}\mtx{S}_{22}].
\end{equation}
Equation (\ref{eq:CPQSidentity}) tells us that the $m\times k$ matrix $\mtx{Q}_{1}\mtx{S}_{11}$
consists precisely of first $k$ columns of $\mtx{A}\mtx{P}$. These columns were the first $k$
columns that were chosen as ``pivots'' in the QR-factorization procedure. They typically form
a good (approximate) basis for the column space of $\mtx{A}$. We consequently define our
$m\times k$ matrix $\mtx{C}$ as this matrix holding the first $k$ pivot columns. Letting $J$
denote the permutation vector associated with the permutation matrix $\mtx{P}$, so that
$$
\mtx{A}\mtx{P} = \mtx{A}(:,J),
$$
we define $J_{\rm s} = J(1:k)$ as the index vector identifying the first $k$ pivots, and set
\begin{equation}
\label{eq:defC}
\mtx{C} = \mtx{A}(:,J_{\rm s}) = \mtx{Q}_{1}\mtx{S}_{11}.
\end{equation}
Now let us rewrite (\ref{eq:CPQSidentity}) by extracting the product $\mtx{Q}_{1}\mtx{S}_{11}$
\begin{equation}
\label{eq:APtmp}
\mtx{A}\mtx{P} =
\mtx{Q}_{1}[\mtx{S}_{11}\ \ |\ \ \mtx{S}_{12}] +
\mtx{Q}_{2}[\mtx{0}     \ \ |\ \ \mtx{S}_{22}] =
\mtx{Q}_{1}\mtx{S}_{11}[\mtx{I}_{k}\ \ |\ \ \mtx{S}_{11}^{-1}\mtx{S}_{12}] +
\mtx{Q}_{2}[\mtx{0}     \ \ |\ \ \mtx{S}_{22}].
\end{equation}
(Remark \ref{remark:Tstability} discusses why $\mtx{S}_{11}$ must be invertible.)
Now define
\begin{equation}
\label{eq:defTZ}
\mtx{T} = \mtx{S}_{11}^{-1}\mtx{S}_{12},
\qquad\mbox{and}\qquad
\mtx{Z} = [\mtx{I}_{k}\ \ |\ \ \mtx{T}]\mtx{P}^{*}
\end{equation}
so that (\ref{eq:APtmp}) can be rewritten (upon right-multiplication by $\mtx{P}^{*}$, which equals $\mtx{P}^{-1}$
since $\mtx{P}$ is unitary) as
\begin{equation}
\label{eq:CZwitherror}
\mtx{A} =
\mtx{C}[\mtx{I}_{k}\ \ |\ \ \mtx{T}]\mtx{P}^{*} +
\mtx{Q}_{2}[\mtx{0}     \ \ |\ \ \mtx{S}_{22}]\mtx{P}^{*} =
\mtx{C}\mtx{Z} +
\mtx{Q}_{2}[\mtx{0}     \ \ |\ \ \mtx{S}_{22}]\mtx{P}^{*}.
\end{equation}
Equation (\ref{eq:CZwitherror}) is precisely the column ID we sought,
with the additional bonus that the remainder term is explicitly identified.
Observe that when the spectral or Frobenius norms are used, the error term
is of \textit{exactly} the same size as the error term obtained from a truncated QR factorization:
$$
\|\mtx{A} - \mtx{C}\mtx{Z}\| =
\|\mtx{Q}_{2}[\mtx{0}     \ \ |\ \ \mtx{S}_{22}]\mtx{P}^{*}\| =
\|\mtx{S}_{22}\| =
\|\mtx{A} - \mtx{Q}_{1}[\mtx{S}_{11}\ \ |\ \ \mtx{S}_{12}]\mtx{P}^{*}\|.
$$

\begin{remark}[Conditioning]
\label{remark:Tstability}
Equation (\ref{eq:APtmp}) involves the quantity $\mtx{S}_{11}^{-1}$ which
prompts the question of whether $\mtx{S}_{11}$ is necessarily invertible,
and what its condition number might be. It is easy to show that whenever
the rank of $\mtx{A}$ is at least $k$, the CPQR algorithm is guaranteed to
result in a matrix $\mtx{S}_{11}$ that is non-singular. (If the rank of
$\mtx{A}$ is $j$, where $j < k$, then the QR factorization process can
detect this and halt the factorization after $j$ steps.) Unfortunately,
$\mtx{S}_{11}$ is typically quite ill-conditioned. The saving grace is
that even though one should expect  $\mtx{S}_{11}$ to be poorly conditioned,
it is often the case that the linear system
\begin{equation}
\label{eq:T_eqn}
\mtx{S}_{11}\mtx{T} = \mtx{S}_{12}
\end{equation}
still has a solution
$\mtx{T}$ whose entries are of moderate size. Informally, one could say
that the directions where $\mtx{S}_{11}$ and $\mtx{S}_{12}$ are small
``line up.'' For standard column pivoted QR, the system (\ref{eq:T_eqn})
will in practice be observed to almost always have a solution $\mtx{T}$
of small size \cite{2005_martinsson_skel}, but counter-examples can be constructed. More sophisticated
pivot selection procedures have been proposed that are \textit{guaranteed}
to result in matrices $\mtx{S}_{11}$ and $\mtx{S}_{12}$ such that
(\ref{eq:T_eqn}) has a good solution; but these are harder to code
and take longer to execute \cite{gu1996}.
\end{remark}

Of course, the row ID can be computed via an entirely analogous process
that starts with a CPQR of the \textit{transpose} of $\mtx{A}$. In other
words, we execute a pivoted Gram-Schmidt orthonormalization process on
the rows of $\mtx{A}$.

Finally, to obtain the double-sided ID, we start with using the CPQR-based
process to build the column ID (\ref{eq:defID1}). Then compute the row ID
by performing Gram-Schmidt on the rows of the tall thin matrix $\mtx{C}$.

The three deterministic algorithms described for computing the three flavors
of ID are summarized in Figure \ref{fig:deterministicID}.

\begin{figure}
\begin{center}
\fbox{\begin{minipage}{100mm}
\textit{Compute a column ID so that $\mtx{A} \approx \mtx{A}(:,J_{\rm s})\,\mtx{Z}$.}\\
\textbf{function} $[J_{\rm s},\ \mtx{Z}] = \texttt{ID\_col}(\mtx{A},k)$\\
\mbox{}\quad$[\mtx{Q},\mtx{S},J] = \texttt{qr}(\mtx{A},0)$;\\
\mbox{}\quad$\mtx{T} = (\mtx{S}(1:k,1:k))^{-1}\mtx{S}(1:k,(k+1):n)$;\\
\mbox{}\quad$\mtx{Z} = \texttt{zeros}(k,n)$\\
\mbox{}\quad$\mtx{Z}(:,J) = [\mtx{I}_{k}\ \mtx{T}]$;\\
\mbox{}\quad$J_{\rm s} = J(1:k)$;

\vspace{2mm}

\textit{Compute a row ID so that $\mtx{A} \approx \mtx{X}\,\mtx{A}(I_{\rm s},:)$.}\\
\textbf{function} $[I_{\rm s},\ \mtx{X}] = \texttt{ID\_row}(\mtx{A},k)$\\
\mbox{}\quad$[\mtx{Q},\mtx{S},\,J] = \texttt{qr}(\mtx{A}^{*},0)$;\\
\mbox{}\quad$\mtx{T} = (\mtx{S}(1:k,1:k))^{-1}\mtx{S}(1:k,(k+1):m)$;\\
\mbox{}\quad$\mtx{X} = \texttt{zeros}(m,k)$\\
\mbox{}\quad$\mtx{X}(J,:) = [\mtx{I}_{k}\ \mtx{T}]^{*}$;\\
\mbox{}\quad$I_{\rm s} = J(1:k)$;

\vspace{2mm}

\textit{Compute a double-sided ID so that $\mtx{A} \approx \mtx{X}\,\mtx{A}(I_{\rm s},J_{\rm s})\,\mtx{Z}$.}\\
\textbf{function} $[I_{\rm s},J_{\rm s},\mtx{X},\mtx{Z}] = \texttt{ID\_double}(\mtx{A},k)$\\
\mbox{}\quad$[J_{\rm s},\ \mtx{Z}] = \texttt{ID\_col}(\mtx{A},k)$;\\
\mbox{}\quad$[I_{\rm s},\ \mtx{X}] = \texttt{ID\_row}(\mtx{A}(:,J_{\rm s}),k)$;

\end{minipage}}
\end{center}
\caption{Deterministic algorithms for computing the column, row, and double-sided ID via the
column pivoted QR factorization. The input is in every case an $m\times n$ matrix $\mtx{A}$
and a target rank $k$. Since the algorithms are based on the CPQR, it is elementary to modify
them to the situation where a tolerance rather than a rank is given. (Recall that the errors
resulting from these ID algorithms are \textit{identical} to the error in the first CPQR factorization
executed.)}
\label{fig:deterministicID}
\end{figure}

\begin{remark}[Partial factorization]
The algorithms for computing interpolatory decompositions shown in
Figure \ref{fig:deterministicID} are wasteful when $k \ll \min(m,n)$ since they
involve a \textit{full} QR factorization, which has complexity $O(mn\min(m,n))$.
This problem is very easily remedied by replacing the full QR factorization by
a \textit{partial} QR factorization, which has cost $O(mnk)$. Such a partial
factorization could take as input either a preset rank $k$, or a tolerance
$\varepsilon$. In the latter case, the factorization would stop once the residual error
$\|\mtx{A} - \mtx{Q}(:,1:k)\mtx{R}(1:k,:)\| = \|\mtx{S}_{22}\| \leq \varepsilon$.
When the QR factorization is interrupted after $k$ steps, the output would still
be a factorization of the form (\ref{eq:defCPQS}), but in this case,
$\mtx{S}_{22}$ would not be upper triangular. This is immaterial since $\mtx{S}_{22}$
is never used. To further accelerate the computation, one can advantageously
use a \textit{randomized CPQR} algorithm, cf.~Sections \ref{sec:randQR} and \ref{sec:randUTV}
or \cite{2015_blockQR_SISC,2015_martinsson_blocked}.
\end{remark}

\subsection{Randomized techniques for computing the ID}
\label{sec:randID}
The ID is particularly well suited to being computed via randomized algorithms.
To describe the ideas, suppose temporarily that $\mtx{A}$ is an $m\times n$
matrix of \textit{exact} rank $k$, and that we have by some means computed
an approximate rank-$k$ factorization
\begin{equation}
\label{eq:YF}
\begin{array}{cccccccccc}
\mtx{A} &=& \mtx{Y} & \mtx{F}. \\
m\times n && m\times k & k\times n
\end{array}
\end{equation}
Once the factorization (\ref{eq:YF}) is available, let us use the algorithm
\texttt{ID\_row} described in Figure \ref{fig:deterministicID} to compute
a row ID $[I_{\rm s},\mtx{X}] = \texttt{ID\_row}(\mtx{Y},k)$ of $\mtx{Y}$ so that
\begin{equation}
\label{eq:Yid}
\begin{array}{cccccccccc}
\mtx{Y} &=& \mtx{X} & \mtx{Y}(I_{\rm s},:). \\
m\times k && m\times k & k\times k
\end{array}
\end{equation}
It then turns out that $\{I_{\rm s},\mtx{X}\}$ is automatically (!) a row ID of $\mtx{A}$ as well.
To see this, simply note that
\begin{multline*}
\mtx{X}\mtx{A}(I_{\rm s},:) = \mbox{\{Use (\ref{eq:YF}) restricted to the rows in $I_{\rm s}$.\}} = \\
\mtx{X}\mtx{Y}(I_{\rm s},:)\mtx{F} = \mbox{\{Use (\ref{eq:Yid}).\}} =
\mtx{Y}\mtx{F} = \mbox{\{Use (\ref{eq:YF}).\}} =
\mtx{A}.
\end{multline*}
The key insight here is very simple, but powerful, so let us spell it out explicitly:
\begin{center}
\fbox{\begin{minipage}{0.95\textwidth}
\textit{Observation:} In order to compute a row ID of a matrix $\mtx{A}$, the only
information needed is a matrix $\mtx{Y}$ whose columns span the column space
of $\mtx{A}$.
\end{minipage}}
\end{center}
As we have seen, the task of finding a matrix $\mtx{Y}$ whose columns form a good
basis for the column space of a matrix is ideally suited to randomized sampling. To
be precise, we showed in Section \ref{sec:stageA} that given a matrix $\mtx{A}$, we
can find a matrix $\mtx{Y}$ whose columns approximately span the column space of
$\mtx{A}$ via the formula $\mtx{Y} = \mtx{A}\mtx{G}$, where $\mtx{G}$ is a tall thin
Gaussian random matrix. The algorithm that results from combining these two insights
is summarized in Figure \ref{fig:randomizedIDpower}.

\begin{figure}
\fbox{\begin{minipage}{0.95\textwidth}
\begin{center}
\textsc{Algorithm: Randomized ID}
\end{center}

\vspace{2mm}

\textit{Inputs:} An $m\times n$ matrix $\mtx{A}$, a target rank $k$, an over-sampling parameter $p$ (say $p=10$),
and a small integer $q$ denoting the number of power iterations taken.

\vspace{2mm}

\textit{Outputs:} An $m\times k$ interpolation matrix $\mtx{X}$ and an index vector $I_{\rm s} \in \mathbb{N}^{k}$
such that $\mtx{A} \approx \mtx{X}\mtx{A}(I_{\rm s},\colon)$.

\vspace{2mm}

\begin{tabbing}
\hspace{8mm} \= \hspace{5mm} \= \hspace{5mm} \= \hspace{15mm}\kill
(1) \> $\mtx{G} = \texttt{randn}(n,k+p)$;\\[1mm]
(2) \> $\mtx{Y} = \mtx{A}\mtx{G}$;\\[1mm]
(3) \> \textbf{for} $j = 1:q$\\[1mm]
(4) \> \> $\mtx{Y}' = \mtx{A}^{*}\mtx{Y};$\\[1mm]
(5) \> \> $\mtx{Y} = \mtx{A}\mtx{Y}';$\\[1mm]
(6) \> \textbf{end for}\\[1mm]
(7) \> Form an ID of the $n\times (k+p)$ sample matrix: $[I_{\rm s},\mtx{X}] = \texttt{ID\_row}(\mtx{Y},k)$.
\end{tabbing}
\end{minipage}}
\caption{An $O(mnk)$ algorithm for computing an interpolative decomposition of $\mtx{A}$ via randomized sampling.
For $q=0$, the scheme is fast
and accurate for matrices whose singular values decay rapidly.
For matrices whose singular values decay slowly, one should pick a larger $q$ (say $q = 1$ or $2$)
to improve accuracy at the cost of longer execution time.
If accuracy better than $\epsilon_{\rm mach}^{1/(2q+1)}$ is desired, then the scheme should be
modified to incorporate orthonormalization as described
in Remark \ref{remark:RSVDP_stabilized}.}
\label{fig:randomizedIDpower}
\end{figure}

The randomized algorithm for computing a row ID shown in Figure \ref{fig:randomizedIDpower}
has complexity $O(mnk)$. We can reduce this complexity to $O(mn\log k)$ by using a structured
random matrix instead of a Gaussian, cf.~Section \ref{sec:SRFT}. The resulting algorithm
is summarized in Figure \ref{fig:randomizedID}.

\begin{figure}
\fbox{\begin{minipage}{0.95\textwidth}
\begin{center}
\textsc{Algorithm: Fast randomized ID}
\end{center}

\vspace{2mm}

\textit{Inputs:} An $m\times n$ matrix $\mtx{A}$, a target rank $k$, and an over-sampling parameter $p$ (say $p=k$).

\vspace{2mm}

\textit{Outputs:} An $m\times k$ interpolation matrix $\mtx{X}$ and an index vector $I_{\rm s} \in \mathbb{N}^{k}$
such that $\mtx{A} \approx \mtx{X}\mtx{A}(I_{\rm s},\colon)$.

\vspace{2mm}

\begin{enumerate}
\item Form an $n\times (k+p)$ SRFT $\mtx{\Omega}$.
\item Form the sample matrix $\mtx{Y} = \mtx{A}\,\mtx{\Omega}$.
\item Form an ID of the $n\times (k+p)$ sample matrix: $[I_{\rm s},\mtx{X}] = \texttt{ID\_row}(\mtx{Y},k)$.
\end{enumerate}
\end{minipage}}
\caption{An $O(mn\log k)$ algorithm for computing an interpolative decomposition of $\mtx{A}$.}
\label{fig:randomizedID}
\end{figure}

\section{Randomized algorithms for computing the CUR decomposition}
\label{sec:CUR}

\subsection{The CUR decomposition}
The so called \textit{CUR-factorization} \cite{2009_mahoney_CUR} is a
``structure preserving'' factorization that is similar to the Interpolative
Decomposition described in Section \ref{sec:ID}.
The CUR factorization approximates  an $m\times n$ matrix $\mtx{A}$ as a product
\begin{equation}
\label{eq:CUR}
\begin{array}{ccccc}
\mtx{A} &\approx & \mtx{C} & \mtx{U} & \mtx{R},\\
m\times n && m\times k & k\times k & k\times n
\end{array}
\end{equation}
where $\mtx{C}$ contains a subset of the columns of $\mtx{A}$ and $\mtx{R}$ contains
a subset of the rows of $\mtx{A}$.
Like the ID, the CUR decomposition offers the ability to preserve properties like sparsity
or non-negativity in the factors of the decomposition, the prospect to
reduce memory requirements, and excellent tools for data interpretation.

The CUR decomposition is often obtained in three steps \cite{2013_mitrovic_cur,2009_mahoney_CUR}:
(1) Some scheme is used to assign a
weight or the so called leverage score (of importance) \cite{1978_hoaglin_leverage_score}
to each column and row in the matrix.
This is typically done either using the $\ell_2$ norms of the columns and rows or
by using the leading singular vectors of $\mtx{A}$ \cite{2008_drineas_relative_error_CUR,2013_shusen_CUR}.
(2) The matrices $\mtx{C}$ and $\mtx{R}$ are constructed via a randomized sampling procedure,
using the leverage scores to assign a sampling probability to each column and row.
(3) The $\mtx{U}$ matrix is computed via:
\begin{equation}
\label{eq:CURformula}
\mtx{U} \approx \mtx{C}^{\dagger} \mtx{A} \mtx{R}^{\dagger},
\end{equation}
with $\mtx{C}^{\dagger}$ and $\mtx{R}^{\dagger}$ being the pseudoinverses of $\mtx{C}$ and $\mtx{R}$.
Non-randomized approaches to computing the CUR decomposition are discussed in
\cite{2016_sorensen_CUR,2014_martinsson_CUR}.

\begin{remark}[Conditioning of CUR]
\label{remark:condCUR}
For matrices whose singular values experience substantial decay, the
accuracy of the CUR factorization can deteriorate due to effects of
ill-conditioning. To simplify slightly, one would normally expect the
leading $k$ singular values of $\mtx{C}$ and $\mtx{R}$ to be
rough approximations to the leading $k$ singular values of $\mtx{A}$,
so that the condition numbers of $\mtx{C}$ and $\mtx{R}$ would be
roughly $\sigma_{1}(\mtx{A})/\sigma_{k}(\mtx{A})$. Since low-rank
factorizations are most useful when applied to matrices whose singular
values decay reasonably rapidly, we would \textit{typically} expect the
ratio  $\sigma_{1}(\mtx{A})/\sigma_{k}(\mtx{A})$ to be large, which is
to say that $\mtx{C}$ and $\mtx{R}$ would be ill-conditioned. Hence, in
the typical case, evaluation of the formula (\ref{eq:CURformula}) can be
expected to result in substantial loss of accuracy due to accumulation
of round-off errors. Observe that the ID does not suffer from this problem;
in (\ref{eq:defID3}), the matrix $\mtx{A}_{\rm skel}$ tends to be
ill-conditioned, but it does not need to be inverted. (The matrices
$\mtx{X}$ and $\mtx{Z}$ are well-conditioned.)
\end{remark}

\subsection{Converting a double-sided ID to a CUR decomposition}
\label{sec:CURheuristics}
We will next describe an algorithm that converts a double-sided ID
to a CUR decomposition. To this end, we assume that the factorization (\ref{eq:defID3})
has been computed using the procedures described in Section \ref{sec:ID}
(either the deterministic or the randomized ones).
In other words, we assume that the index vectors $I_{\rm s}$ and $J_{\rm s}$,
and the basis matrices $\mtx{X}$ and $\mtx{Z}$, are all available.
We then define $\mtx{C}$ and $\mtx{R}$ in the natural way as
\begin{equation}
\label{eq:choice_C_and_R}
\mtx{C} = \mtx{A}(:,J_{\rm s})
\quad\text{and}\quad
\mtx{R} = \mtx{A}(I_{\rm s},:).
\end{equation}
Consequently, $\mtx{C}$ and $\mtx{R}$ are respectively subsets of columns and of
rows of $\mtx{A}$. The index vectors $I_{\rm s}$
and $J_{\rm s}$ are determined by the column pivoted QR factorizations, possibly
combined with a randomized projection step for computational efficiency.
It remains to construct a $k\times k$ matrix $\mtx{U}$ such that
\begin{equation}
\label{eq:netflix2}
\mtx{A} \approx \mtx{C}\,\mtx{U}\,\mtx{R}.
\end{equation}
Now recall that, cf.~(\ref{eq:defID1}),
\begin{equation}
\label{eq:netflix1}
\mtx{A} \approx \mtx{C}\,\mtx{Z}.
\end{equation}
By inspecting (\ref{eq:netflix1}) and (\ref{eq:netflix2}), we find that we
achieve our objective if we determine a matrix $\mtx{U}$ such that
\begin{equation}
\label{eq:overd}
\begin{array}{cccccccccccccc}
\mtx{U} & \mtx{R} &=& \mtx{Z}.\\
k\times k & k\times n && k\times n
\end{array}
\end{equation}
Unfortunately, (\ref{eq:overd}) is an over-determined system, but at least intuitively,
it seems plausible that it should have a fairly accurate solution, given that the
rows of $\mtx{R}$ and the rows of $\mtx{Z}$ should, by construction, span roughly the
same space (namely, the space spanned by the $k$ leading right singular vectors
of $\mtx{A}$). Solving (\ref{eq:overd}) in the least-square sense, we arrive at our
definition of $\mtx{U}$:
\begin{equation}
\label{eq:defU}
\mtx{U} := \mtx{Z}\mtx{R}^{\dagger}.
\end{equation}

The resulting algorithm is summarized in Figure \ref{fig:randomizedCUR}.

\begin{figure}
\fbox{\begin{minipage}{0.95\textwidth}
\begin{center}
\textsc{Algorithm: Randomized CUR}
\end{center}

\vspace{2mm}

\textit{Inputs:} An $m\times n$ matrix $\mtx{A}$, a target rank $k$, an over-sampling parameter $p$ (say $p=10$),
and a small integer $q$ denoting the number of power iterations taken.

\vspace{2mm}

\textit{Outputs:} A $k\times k$ matrix $\mtx{U}$ and index vectors $I_{\rm s}$ and $J_{\rm s}$
of length $k$ such that $\mtx{A} \approx \mtx{A}(\colon,J_{\rm s})\,\mtx{U}\,\mtx{A}(I_{\rm s},\colon)$.

\vspace{2mm}

\begin{tabbing}
\hspace{8mm} \= \hspace{5mm} \= \hspace{5mm} \= \hspace{15mm}\kill
(1) \> $\mtx{G} = \texttt{randn}(k+p,m)$;\\[1mm]
(2) \> $\mtx{Y} = \mtx{G}\mtx{A}$;\\[1mm]
(3) \> \textbf{for} $j = 1:q$\\[1mm]
(4) \> \> $\mtx{Z} = \mtx{Y}\mtx{A}^{*};$\\[1mm]
(5) \> \> $\mtx{Y} = \mtx{Z}\mtx{A};$\\[1mm]
(6) \> \textbf{end for}\\[1mm]
(7) \> Form an ID of the $(k+p)\times n$ sample matrix: $[I_{\rm s},\mtx{Z}] = \texttt{ID\_col}(\mtx{Y},k)$.\\[1mm]
(8) \> Form an ID of the $m\times k$ matrix of chosen columns: $[I_{\rm s},\sim] =
                                                                \texttt{ID\_row}(\mtx{A}(:,J_{\rm s}),k)$.\\[1mm]
(9) \> Determine the matrix $\mtx{U}$ by solving the least squares equation
       $\mtx{U}\mtx{A}(I_{\rm s},:) = \mtx{Z}$.
\end{tabbing}
\end{minipage}}
\caption{A randomized algorithm for computing a CUR decomposition of $\mtx{A}$ via randomized sampling.
For $q=0$, the scheme is fast and accurate for matrices whose singular values decay rapidly.
For matrices whose singular values decay slowly, one should pick a larger $q$ (say $q = 1$ or $2$)
to improve accuracy at the cost of longer execution time.
If accuracy better than $\epsilon_{\rm mach}^{1/(2q+1)}$ is desired, then the scheme should be
modified to incorporate orthonormalization as described
in Remark \ref{remark:RSVDP_stabilized}.}
\label{fig:randomizedCUR}
\end{figure}

\section{Adaptive rank determination with updating of the matrix}
\label{sec:acc_given}

\subsection{Problem formulation}
Up to this point, we have assumed that the rank $k$ is given as an input
variable to the factorization algorithm. In practical usage, it is common
that we are given instead a matrix $\mtx{A}$ and a computational tolerance
$\varepsilon$, and our task is then to determine a matrix $\mtx{A}_{k}$ of
rank $k$ such that $\|\mtx{A} - \mtx{A}_{k}\| \leq \varepsilon$.

The techniques described in this section are designed for dense matrices
stored in RAM. They directly update the matrix, and come with a firm
guarantee that the computed low rank approximation is within distance
$\varepsilon$ of the original matrix. There are many situations where
direct updating is not feasible and we can in practice only interact
with the matrix via the matrix-vector multiplication (e.g., very large
matrices stored out-of-core, sparse matrices, matrices that are defined
implicitly). Section \ref{sec:rand_noupdating} describes algorithms designed
for this environment that use randomized sampling techniques to
\textit{estimate} the approximation error.

Recall that for the case where a computational tolerance is given (rather
than a rank), the optimal solution is given by the SVD. Specifically,
let $\{\sigma_{j}\}_{j=1}^{\min(m,n)}$ be the singular values of $\mtx{A}$,
and let $\varepsilon$ be a given tolerance.
Then the minimal rank $k$ for which there exists a matrix $\mtx{B}$ of rank $k$
that is within distance $\varepsilon$ of $\mtx{A}$, is the
is the smallest integer $k$ such that $\sigma_{k+1} \leq \varepsilon$.
The algorithms described here will determine a $k$ that is not necessarily optimal, but
is typically fairly close.

\subsection{A greedy updating algorithm}
\label{sec:greedy_basic}
Let us start by describing a general algorithmic template for how
to compute an approximate rank-$k$ approximate factorization of a matrix.
To be precise, suppose that we are given an $m\times n$ matrix $\mtx{A}$,
and a computational tolerance $\varepsilon$. Our objective is then to
determine an integer $k \in \{1,2,\dots,\min(m,n)\}$, an $m\times k$
orthonormal matrix $\mtx{Q}_{k}$, and a $k\times n$ matrix $\mtx{B}_{k}$ such that
\[
\begin{array}{ccccccccccccccc}
\|&\mtx{A} &-& \mtx{Q}_{k} & \mtx{B}_{k} & \| \leq \varepsilon.\\
& m\times n && m\times k & k\times n
\end{array}
\]
Figure \ref{fig:greedy_template} outlines a template for how one might in a greedy fashion
build the matrices $\mtx{Q}_{k}$ and $\mtx{B}_{k}$ by adding one column to $\mtx{Q}_{k}$,
and one row to $\mtx{B}_{k}$, at each step.

\begin{figure}
\begin{center}
\fbox{\begin{minipage}{70mm}
\begin{tabbing}
\hspace{8mm} \= \hspace{5mm} \= \hspace{5mm} \= \hspace{15mm}\kill
(1)  \> $\mtx{Q}_{0} = [\ ]$; $\mtx{B}_{0} = [\ ]$; $\mtx{A}_{0} = \mtx{A}$; $k=0$;\\[1mm]
(2)  \> \textbf{while} $\|\mtx{A}_{k}\| > \varepsilon$\\[1mm]
(3)  \> \> $k = k+1$\\[1mm]
(4)  \> \> Pick a vector $\vct{y} \in \mbox{Ran}(\mtx{A}_{k-1})$\\[1mm]
(5)  \> \> $\vct{q} = \vct{y} / \|\vct{y}\|;$\\[1mm]
(6)  \> \> $\vct{b} = \vct{q}^{*}\mtx{A}_{k-1};$\\[1mm]
(7)  \> \> $\mtx{Q}_{k} = [\mtx{Q}_{k-1}\ \vct{q}];$\\[1mm]
(8)  \> \> $\mtx{B}_{k} = \left[\begin{array}{c} \mtx{B}_{k-1} \\ \vct{b} \end{array}\right];$\\[1mm]
(9)  \> \> $\mtx{A}_{k} = \mtx{A}_{k-1} - \vct{q}\vct{b};$\\[1mm]
(10) \> \textbf{end for}
\end{tabbing}
\end{minipage}}
\end{center}
\caption{A greedy algorithm for building a low-rank approximation to a given $m\times n$
matrix $\mtx{A}$ that is accurate to within a given precision $\varepsilon$. To be precise,
the algorithm determines an integer $k$, an $m\times k$ orthonormal matrix $\mtx{Q}_{k}$ and a
$k\times n$ matrix $\mtx{B}_{k} = \mtx{Q}_{k}^{*}\mtx{A}$ such that
$\|\mtx{A} - \mtx{Q}_{k}\mtx{B}_{k}\| \leq \varepsilon$. One can easily verify that
after step $j$, we have $\mtx{A} = \mtx{Q}_{j}\mtx{B}_{j} + \mtx{A}_{j}$.}
\label{fig:greedy_template}
\end{figure}

The algorithm described in Figure \ref{fig:greedy_template} is a generalization of
the classical Gram-Schmidt procedure.
The key to understanding how the algorithm works
is provided by the identity
\[
\mtx{A} = \mtx{Q}_{j}\mtx{B}_{j} + \mtx{A}_{j},\qquad j = 0,1,2,\dots,k.
\]
The computational efficiency and accuracy of the algorithm depend crucially on the
strategy for picking the vector $\vct{y}$ on line (4). Let us consider three possibilities:

\subsubsection*{Pick the largest remaining column} Suppose we instantiate line (4) by
letting $\vct{y}$ be simply the largest column of the remainder matrix $\mtx{A}_{k-1}$.
\begin{center}
\begin{minipage}{70mm}
\begin{tabbing}
\hspace{8mm} \= \hspace{5mm} \= \hspace{5mm} \= \hspace{15mm}\kill
(4)  \> \> Set $j_{k} = \mbox{argmax}\{\|\mtx{A}_{k-1}(:,j)\|\,\colon\,j = 1,2,\dots,n\}$
           and then $\vct{y} = \mtx{A}_{k-1}(:,j_{k})$.
\end{tabbing}
\end{minipage}
\end{center}
With this choice, the algorithm in Figure \ref{fig:greedy_template} is \textit{precisely}
column pivoted Gram-Schmidt (CPQR). This algorithm is reasonably efficient, and often leads to
fairly close to optimal low-rank approximation. For instance, when the singular values
of $\mtx{A}$ decay rapidly, CPQR determines a numerical rank $k$ that is
typically reasonably close to the theoretically exact $\varepsilon$-rank. However, this is
not always the case even when the singular values decay rapidly, and the results can be
quite poor when the singular values decay slowly. (A striking illustration of how suboptimal
CPQR can be for purposes of low-rank approximation is provided by the famous ``Kahan
counter-example,'' see \cite[Sec.~5]{1966_kahan_NLA}.)

\subsubsection*{Pick the locally optimal vector} A choice that is natural and
conceptually simple is to pick the vector $\vct{y}$ by solving the obvious
minimization problem:
\begin{center}
\begin{minipage}{70mm}
\begin{tabbing}
\hspace{8mm} \= \hspace{5mm} \= \hspace{5mm} \= \hspace{15mm}\kill
(4)  \> \> $\vct{y} = \mbox{argmin}\{\|\mtx{A}_{k-1} - \vct{y}\vct{y}^{*}\mtx{A}_{k-1}\|\,\colon\,
\|\vct{y}\| = 1\}$.
\end{tabbing}
\end{minipage}
\end{center}
With this choice, the algorithm will produce matrices that attain the theoretically
optimal precision
\[
\|\mtx{A} - \mtx{Q}_{j}\mtx{B}_{j}\| = \sigma_{j+1}.
\]
This tells us that the greediness of the algorithm is not a problem.
However, this strategy is impractical since solving the local minimization
problem is computationally hard.

\subsubsection*{A randomized selection strategy} Suppose now that we pick $\vct{y}$ by
forming a linear combination of the columns of $\mtx{A}_{k-1}$ with the expansion weights
drawn from a normalized Gaussian distribution:
\begin{center}
\begin{minipage}{70mm}
\begin{tabbing}
\hspace{8mm} \= \hspace{5mm} \= \hspace{5mm} \= \hspace{15mm}\kill
(4)  \> \> Draw a Gaussian random vector $\vct{g} \in \mathbb{R}^{n}$ and set
$\vct{y} = \mtx{A}_{k-1}\vct{g}$.
\end{tabbing}
\end{minipage}
\end{center}
With this choice, the algorithm becomes logically equivalent to the basic randomized SVD
given in Figure \ref{fig:RSVD}. This means that this choice often leads to a factorization
that is close to optimally accurate, and is also computationally efficient. One can attain
higher accuracy by trading away some computational efficiency and incorporate a couple of
steps of power iteration, and choosing
$\vct{y} = \bigl(\mtx{A}_{k-1}\mtx{A}_{k-1}^{*}\bigr)^{q}\mtx{A}_{k-1}\vct{g}$
for some small integer $q$, say $q=1$ or $q=2$.

\subsection{A blocked updating algorithm}
A key benefit of the randomized greedy algorithm described in Section
\ref{sec:greedy_basic} is that it can easily be \textit{blocked}. In
other words, given a block size $b$, we can at each step of the
iteration draw a set of $b$ Gaussian random vectors, compute the
corresponding sample vectors, and then extend the factors $\mtx{Q}$
and $\mtx{B}$ by adding $b$ columns and $b$ rows at a time,
respectively. The resulting algorithm is shown in Figure
\ref{fig:greedy_blocked}.

\begin{figure}
\begin{center}
\fbox{
\begin{minipage}{120mm}
\begin{tabbing}
\hspace{15mm} \= \hspace{15mm} \= \hspace{15mm} \= \hspace{15mm}\kill
(1) \> $\mtx{Q} = [\ ]$; $\mtx{B} = [\ ]$; \\[1mm]
(2) \> \textbf{while} $\|\mtx{A}\| > \varepsilon$\\[1mm]
(3) \> \> Draw an $n\times b$ Gaussian random matrix $\mtx{G}$.\\[1mm]
(4) \> \> Compute the $m\times b$ matrix $\mtx{Q}_{\rm new} = \texttt{qr}(\mtx{A}\mtx{G},0)$.\\[1mm]
(5) \> \> $\mtx{B}_{\rm new} = \vct{Q}_{\rm new}^{*}\mtx{A}$\\[1mm]
(6) \> \> $\mtx{Q} = [\mtx{Q}\ \mtx{Q}_{\rm new}]$\\[1mm]
(7) \> \> $\mtx{B} = \left[\begin{array}{c} \mtx{B} \\ \mtx{B}_{\rm new} \end{array}\right]$\\[1mm]
(8) \> \> $\mtx{A} = \mtx{A} - \mtx{Q}_{\rm new}\mtx{B}_{\rm new}$\\[1mm]
(9) \> \textbf{end while}
\end{tabbing}
\end{minipage}}
\end{center}
\caption{A greedy algorithm for building a low-rank approximation to a given $m\times n$
matrix $\mtx{A}$ that is accurate to within a given precision $\varepsilon$. This algorithm
is a blocked analogue of the method described in Figure \ref{fig:greedy_template} and takes as
input a block size $b$. Its output is an orthonormal matrix $\mtx{Q}$ of size $m\times k$ (where
$k$ is a multiple of $b$) and a $k\times n$ matrix $\mtx{B}$ such that
$\|\mtx{A} - \mtx{Q}\mtx{B}\| \leq \varepsilon$. For higher accuracy, one can incorporate
a couple of steps of power iteration and set $\mtx{Q}_{\rm new} = \texttt{qr}((\mtx{A}\mtx{A}^{*})^{q}\mtx{A}\mtx{R},0)$
on Line (4).}
\label{fig:greedy_blocked}
\end{figure}

\subsection{Evaluating the norm of the residual}
The algorithms described in this section contain one step that could be computationally expensive unless some care is exercised.
The potential problem concerns the evaluation of the norm of the remainder matrix
$\|\mtx{A}_{k}\|$ (cf.~Line (2) in Figures \ref{fig:greedy_template}) at each step of the iteration. When the Frobenius
norm is used, this evaluation can be done very efficiently, as follows: When the
computation starts, evaluate the norm of the input matrix
$$
a = \|\mtx{A}\|_{\rm Fro}.
$$
Then observe that after step $j$ completes, we have
$$
\mtx{A} =
\underbrace{\mtx{Q}_{j}\mtx{B}_{j}}_{=\mtx{Q}_{j}\mtx{Q}_{j}^{*}\mtx{A}} +
\underbrace{\mtx{A}_{j}}_{=(\mtx{I} - \mtx{Q}_{j}\mtx{Q}_{j}^{*}\bigr)\mtx{A}}.
$$
Since the columns in the first term all lie in $\mbox{Col}(\mtx{Q}_{j})$, and the columns of
the second term all lie in $\mbox{Col}(\mtx{Q}_{j})^{\perp}$, we now find that
$$
\|\mtx{A}\|_{\rm Fro}^{2} =
\|\mtx{Q}_{j}\mtx{B}_{j}\|_{\rm Fro}^{2} + \|\mtx{A}_{j}\|_{\rm Fro}^{2} =
\|\mtx{B}_{j}\|_{\rm Fro}^{2} + \|\mtx{A}_{j}\|_{\rm Fro}^{2},
$$
where in the last step we used that $\|\mtx{Q}_{j}\mtx{B}_{j}\|_{\rm Fro} = \|\mtx{B}_{j}\|_{\rm Fro}$
since $\mtx{Q}_{j}$ is orthonormal.
In other words, we can easily compute $\|\mtx{A}_{j}\|_{\rm Fro}$ via the identity
$$
\|\mtx{A}_{j}\|_{\rm Fro} = \sqrt{a^{2} - \|\mtx{B}_{j}\|_{\rm Fro}^{2}}.
$$
The idea here is related to ``down-dating'' schemes for computing column norms when
executing a column pivoted QR factorization, as described in, e.g.,
\cite[Chapter 5, Section 2.1]{1998_stewart_volume1}.

When the spectral norm is used, one could use a power iteration to compute an estimate of
the norm of the matrix. Alternatively, one can use the randomized procedure described in
Section \ref{sec:rand_noupdating} which is faster, but less accurate.

\section{Adaptive rank determination without updating the matrix}
\label{sec:rand_noupdating}
The techniques described in Section \ref{sec:acc_given} for computing a low
rank approximation to a matrix $\mtx{A}$ that is valid to a \textit{given tolerance} (as opposed to
a \textit{given rank}) are highly computationally efficient whenever the matrix $\mtx{A}$
itself can be easily updated (e.g.~a dense matrix stored in RAM). In this section,
we describe algorithms for solving the ``given tolerance'' problem that do not
need to explicitly update the matrix; this comes in handy for sparse matrices,
for matrices stored out-of-core, for matrices defined implicitly, etc.
In such a situation, it often works well to use randomized methods to
estimate the norm of the residual matrix $\mtx{A} - \mtx{Q}\mtx{B}$. The framing
for such a randomized estimator is that given a tolerated risk probability $p$,
we can cheaply compute a bound for $\|\mtx{A} - \mtx{Q}\mtx{B}\|$ that is valid with probability at
least $1-p$.

\begin{figure}
\begin{center}
\fbox{\begin{minipage}{70mm}
\begin{tabbing}
\hspace{8mm} \= \hspace{5mm} \= \hspace{5mm} \= \hspace{15mm}\kill
(1)  \> $\mtx{Q}_{0} = [\ ]$; $\mtx{B}_{0} = [\ ]$; \\[1mm]
(2)  \> \textbf{for} $j = 1,2,3,\dots$\\[1mm]
(3)  \> \> Draw a Gaussian random vector $\vct{g}_{j} \in \mathbb{R}^{n}$ and set
           $\vct{y}_{j} = \mtx{A}\vct{g}_{j}$\\[1mm]
(4)  \> \> Set $\vct{z}_{j} = \vct{y}_{j}  - \mtx{Q}_{j-1}\mtx{Q}_{j-1}^{*}\vct{y}_{j}$
           and then $\vct{q}_{j} = \vct{z}_{j}/|\vct{z}_{j}|$.\\[1mm]
(5)  \> \> $\mtx{Q}_{j} = [\mtx{Q}_{j-1}\ \vct{q}_{j}]$\\[1mm]
(6)  \> \> $\mtx{B}_{j} = \left[\begin{array}{c} \mtx{B}_{j-1} \\ \vct{q}_{j}^{*}\mtx{A} \end{array}\right]$\\[1mm]
(7)  \> \textbf{end for}
\end{tabbing}
\end{minipage}}
\end{center}
\caption{A randomized range finder that builds an orthonormal basis $\{\vct{q}_{1},\vct{q}_{2},\vct{q}_{3},\dots\}$
for the range of $\mtx{A}$ one vector at a time. This algorithm is mathematically equivalent to the basic
RSVD in Figure \ref{fig:RSVD} in the sense that if $\mtx{G} = [\vct{g}_{1}\ \vct{g}_{2}\ \vct{g}_{3}\ \dots]$,
then the vectors $\{\vct{q}_{j}\}_{j=1}^{p}$ form an orthonormal basis for $\mtx{A}\mtx{G}(:,1:p)$ for both methods.
Observe that $\vct{z}_{j} =
\bigl(\mtx{A} - \mtx{Q}_{j-1}\mtx{Q}_{j-1}^{*}\mtx{A}\bigr)\vct{g}_{j}$, cf.~(\ref{eq:rookendgamge}).}
\label{fig:singlevecrangefinder}
\end{figure}

As a preliminary step in deriving the update-free scheme, let us reformulate
the basic RSVD in Figure \ref{fig:RSVD} as the sequential algorithm shown in
Figure \ref{fig:singlevecrangefinder} that builds the matrices
$\mtx{Q}$ and $\mtx{B}$ one vector at a time. We observe that this
method is similar to the greedy template shown in Figure \ref{fig:greedy_template},
except that there is not an immediately obvious way to tell when $\|\mtx{A} - \mtx{Q}_{j}\mtx{B}_{j}\|$
becomes small enough. However, it is possible to \textit{estimate} this quantity quite
easily. The idea is that once $\mtx{Q}_{j}$ becomes large enough to capture ``most'' of
the range of $\mtx{A}$, then the sample vectors $\vct{y}_{j}$ drawn will all approximately
lie in the span of $\mtx{Q}_{j}$, which is to say that the projected vectors
$\vct{z}_{j}$ will become very small. In other words, once we start to see a sequence
of vectors $\vct{z}_{j}$ that are all very small, we can reasonably deduce that the
basis we have on hand very likely covers most of the range of $\mtx{A}$.

To make the discussion in the previous paragraph more mathematically rigorous, let
us first observe that each projected vector $\mtx{z}_{j}$ satisfies the relation
\begin{equation}
\label{eq:rookendgamge}
\vct{z}_{j} =
\mtx{y}_{j} - \mtx{Q}_{j-1}\mtx{Q}_{j-1}^{*}\vct{y}_{j} =
\mtx{A}\vct{g}_{j} - \mtx{Q}_{j-1}\mtx{Q}_{j-1}^{*}\mtx{A}\vct{g}_{j} =
\bigl(\mtx{A} - \mtx{Q}_{j-1}\mtx{Q}_{j-1}^{*}\mtx{A}\bigr)\vct{g}_{j}.
\end{equation}
Next, we use the fact that if $\mtx{T}$ is any matrix, then by looking at the
magnitude of $\|\mtx{T}\vct{g}\|$, where $\vct{g}$ is a Gaussian random vector,
we can deduce information about the spectral norm of $\mtx{T}$.
The precise result that we need is the following,
cf.~\cite{2007_woolfe_liberty_rokhlin_tygert}*{Sec.~3.4} and \cite{2011_martinsson_randomsurvey}*{Lemma 4.1}.

\begin{lemma}
\label{thm:aposteriori}
Let $\mtx{T}$ be a real $m\times n$ matrix.
Fix a positive integer $r$ and a real number $\alpha \in (0,1)$.
Draw an independent family $\{ \vct{g}_{i} : i = 1, 2, \dots, r \}$
of standard Gaussian vectors.  Then
\begin{equation*}
\|\mtx{T}\|
    \leq \frac{1}{\alpha } \sqrt{\frac{2}{\pi}} \max_{i = 1, \dots, r}
    \|\mtx{T}\vct{g}_{i} \|
\end{equation*}
with probability at least $1 - \alpha^{r}$.
\end{lemma}

In applying this result, we set $\alpha = 1/10$, whence it follows that
if $\|\vct{z}_{j}\|$ is smaller than the resulting threshold for $r$ vectors
in a row, then $\|\mtx{A} - \mtx{Q}_{j}\mtx{B}_{j}\| \leq \varepsilon$
with probability at least $1 - 10^{-r}$. The resulting algorithm is shown
in Figure \ref{alg:adaptive2}. Observe that choosing $\alpha=1/10$
will work well only if the singular values decay reasonably fast.

\newcommand{\smnorm}[2]{\|#2\|}

\begin{figure}
\begin{center}
\fbox{
\begin{minipage}{.9\textwidth}
\begin{tabbing}
\hspace{8mm} \= \hspace{8mm} \= \hspace{8mm} \= \hspace{8mm} \= \kill
 (1) \> Draw standard Gaussian vectors $\vct{g}_{1}, \dots, \vct{g}_{r}$ of length $n$.\\
 (2) \> For $i = 1,2,\dots,r$, compute $\vct{y}_{i} = \mtx{A}\vct{g}_{i}$.\\
 (3) \> $j=0$.\\
 (4) \> $\mtx{Q}_{0} = [\ ]$, the $m\times 0$ empty matrix. \\
 (5) \> \textbf{while} $\displaystyle
         \max\left\{\smnorm{}{\vct{y}_{j+1}},\smnorm{}{\vct{y}_{j+2}},\dots,\smnorm{}{\vct{y}_{j+r}} \right\} >
\varepsilon/(10\sqrt{2/\pi})$,\\
 (6) \> \> $j = j + 1$.\\
 (7) \> \> Overwrite $\vct{y}_{j}$ by $\vct{y}_{j} - \mtx{Q}_{j-1}(\mtx{Q}_{j-1})^{*}\vct{y}_{j}$.\\
 (8) \> \> $\vct{q}_{j} = \vct{y}_{j}/|\vct{y}_{j}|$.\\
 (9) \> \> $\mtx{Q}_{j} = [\mtx{Q}_{j-1}\ \vct{q}_{j}]$.\\
 (10) \> \> Draw a standard Gaussian vector $\vct{g}_{j+r}$ of length $n$.\\
 (11) \> \> $\vct{y}_{j+r} = \left(\mtx{I} - \mtx{Q}_{j}(\mtx{Q}_{j})^{*}\right)\mtx{A}\vct{g}_{j+r}$.\\
 (12) \> \> \textbf{for} $i = (j+1),(j+2),\dots,(j+r-1)$,\\
 (13) \> \> \> Overwrite $\vct{y}_{i}$ by $\vct{y}_{i} -
                          \vct{q}_{j}\,\langle \vct{q}_{j},\,\vct{y}_{i}\rangle$.\\
 (14) \> \> \textbf{end for}\\
 (15) \> \textbf{end while}\\
 (16) \> $\mtx{Q} = \mtx{Q}_{j}$.
\end{tabbing}
\end{minipage}}
\end{center}
\caption{A randomized range finder.
Given an $m\times n$ matrix $\mtx{A}$, a tolerance $\varepsilon$,
and an integer $r$, the algorithm computes an orthonormal matrix
$\mtx{Q}$ such that $\|\mtx{A} - \mtx{Q}\mtx{Q}^{*}\mtx{A}\| \leq \varepsilon$ holds with probability at
least $1 - \min\{m,n\} 10^{-r}$. Line (7) is mathematically redundant but improves orthonormality of $\mtx{Q}$
in the presence of round-off errors.
(Adapted from Algorithm 4.2 of \cite{2011_martinsson_randomsurvey}.)}
\label{alg:adaptive2}
\end{figure}

\begin{remark}
While the proof of Lemma \ref{thm:aposteriori} is outside the scope of this tutorial,
it is perhaps instructive to prove a much simpler related result that says that if
$\mtx{T}$ is any $m\times n$ matrix, and $\vct{g}\in \mathbb{R}^{n}$
is a standard Gaussian vector, then
\begin{equation}
\label{eq:infomm}
\mathbb{E}[\|\mtx{T}\vct{g}\|^{2}] = \|\mtx{T}\|^{2}_{\rm Fro}.
\end{equation}
To prove (\ref{eq:infomm}), let $\mtx{T}$ have the singular value decomposition
$\mtx{T} = \mtx{U}\mtx{D}\mtx{V}^{*}$,
and set $\tilde{\vct{g}} = \mtx{V}^{*}\vct{g}$. Then
\[
\|\mtx{T}\vct{g}\|^{2} =
\|\mtx{U}\mtx{D}\mtx{V}^{*}\vct{g}\|^{2} =
\|\mtx{U}\mtx{D}\tilde{\vct{g}}\|^{2} =
\mbox{\{$\mtx{U}$ is unitary\}} =
\|\mtx{D}\tilde{\vct{g}}\|^{2} =
\sum_{j=1}^{n}\sigma_{j}^{2}\tilde{g}_{j}^{2}.
\]
Then observe that since the distribution of Gaussian vectors is rotationally invariant, the
vector $\tilde{\vct{g}} = \mtx{V}^{*}\vct{g}$ is also a standardized Gaussian vector, and
so $\mathbb{E}[\tilde{g}_{j}^{2}] = 1$ for every $j$. Since the variables $\{\tilde{g}_{j}\}_{j=1}^{n}$
are independent, and $\mathbb{E}[\tilde{g}_{j}^{2}]=1$ for every $j$, it follows that
\[
\mathbb{E}[\|\mtx{T}\vct{g}\|^{2}] =
\mathbb{E}\left[\sum_{j=1}^{n}\sigma_{j}^{2}\tilde{g}_{j}^{2}\right] =
\sum_{j=1}^{n}\sigma_{j}^{2}\mathbb{E}[\tilde{g}_{j}^{2}] =
\sum_{j=1}^{n}\sigma_{j}^{2} =
\|\mtx{T}\|_{\rm Fro}^{2}.
\]
\end{remark}

\section{\`A posteriori error bounds, and certificates of accuracy}

The idea of using a randomized method to estimate the norm of a matrix
in Section \ref{sec:rand_noupdating} can be used as a powerful tool
for inexpensively computing a ``certificate of accuracy'' that assures a
user that the output of a matrix factorization algorithm is accurate. This
becomes particularly valuable when ``risky'' randomized algorithms such as those
based on, e.g., fast Johnson-Lindenstrauss transforms (see Section \ref{sec:SRFT}) are used.

To illustrate with a very simple example, suppose we have at our disposal
an algorithm
$$
\mtx{A}_{\rm approx} = \texttt{factorizefast}(\mtx{A},\varepsilon)
$$
that given a matrix $\mtx{A}$ and a tolerance $\varepsilon$ computes a
low rank approximation $\mtx{A}_{\rm approx}$ (which will of course be
returned in terms of the thin factors). Think of this algorithm as typically
doing a good job, but that it is not entirely reliable in the sense that the
likelihood of a substantial error is non-negligible. What we want in this
case is an estimate for the error $\|\mtx{A} - \mtx{A}_{\rm approx}\|$.
Lemma \ref{thm:aposteriori} provides us with one method for doing this
efficiently. Simply draw a \textit{very} thin Gaussian random matrix
$\mtx{G}$ with, say, only 10 columns. Then if we evaluate the matrix
$$
\mtx{E} = \bigl(\mtx{A} - \mtx{A}_{\rm approx}\bigr)\mtx{G},
$$
we will immediately get a highly reliable bound on the error
$\|\mtx{A} - \mtx{A}_{\rm approx}\|$ from Lemma \ref{thm:aposteriori}.
For instance, choosing $\alpha=0.1$, we find that the bound
$$
\|\mtx{A} - \mtx{A}_{\rm approx}\| \leq 10\sqrt{\frac{2}{\pi}}\,\max_{i}\|\mtx{E}(:,i)\|
$$
holds with probability at least $1 - 10^{-10}$.

Let us next discuss how this ``certificate of accuracy'' idea can be deployed
in some typical environment.

\textit{(1) Fast Johnson-Lindenstrauss transforms:} Let us recall the
computational setup of the $O(mn\log k)$ complexity method described in
Section \ref{sec:SRFT}: We use a ``structured'' random matrix
$\mtx{\Omega} \in \mathbb{R}^{n\times \ell}$ that allows us to
rapidly evaluate the sample matrix $\mtx{Y} = \mtx{A}\mtx{\Omega}$,
that is used to form a basis for the column space. We discussed the particular
case of a subsampled randomized Fourier transform, but there are many other
``fast'' random projections to choose from. Methods of this type \textit{typically}
perform almost as well as Gaussian random matrices, but come with far weaker
performance guarantees since one can build an adversarial matrix $\mtx{A}$ for
which the method will fail. The certificate of accuracy idea will in this environment
provide us with a method that has $O(mn\log k)$ complexity, has only very slightly
higher practical execution time as the basic method, and produces an answer that is
in every way just as reliable as the answer resulting from using a basic Gaussian random matrix.
Observe that in the event that an unacceptable large error has been produced, one can
often simply repeat the process using a different instantiation of $\mtx{\Omega}$.

\textit{(2) Single-pass algorithms:}
The techniques described in Section \ref{sec:singlepass} for how to factorize a matrix
that cannot be stored share some characteristics of methods based on Fast Johnson-Lindenstrauss
transforms: The output is in the typical situation accurate enough, but the likelihood of an
uncharacteristically large error is much higher than the traditional ``two-pass'' randomized
SVD with a Gaussian random matrix. In the single-pass environment, the certificate of accuracy
can be implemented by simply building up a secondary sample matrix that records the action of
$\mtx{A}$ on a very thin Gaussian matrix $\mtx{G}_{\rm cert}$ drawn separately for purposes of
authentication. At the end of the process, we can then simply compare the sample matrix
$\mtx{A}\mtx{G}_{\rm cert}$ that we built up over the single pass over the matrix against a
matrix $\mtx{A}_{\rm approx}\mtx{G}_{\rm cert}$ that we can compute using the computed factors
of $\mtx{A}_{\rm approx}$ to provide an estimate for the error. In this case, a user can of course
not go back for a ``second bite of the apple'' in the event that a large error has been produced,
but will at least get a warning that the computed approximation is unreliable.

\textit{(3) Subsampling based on prior knowledge of $\mtx{A}$:}
In the literature on randomized linear algebra, a strand of research that
has attracted much interest is the idea to use a sample matrix $\mtx{Y}$
that consists of $\ell$ randomly chosen columns of $\mtx{A}$ itself to
form a basis for its column space \cite{2011_mahoney_survey,2019_strang_linear_algebra_data}.
For special classes of matrices, it may work well to draw  the samples from
a uniform distribution, but in general, one needs to first do some pre-processing
to obtain sampling probabilities that improve the likelihood of success. A very
simple idea is to use the column norms of $\mtx{A}$ as weights; these are cheap to compute but
generally do not contain enough information. A more sophisticated idea is to use
so called ``leverage scores''  which are the column norms of a matrix whose rows are
the dominant right singular vectors \cite[Sec.~4.3]{2018_PCMI_mahoney_drineas}. These vectors are of course typically not available,
but techniques for estimating the leverage scores have been proposed. Methods of this
type can be very fast, and are for certain classes of matrices highly reliable. In cases
where there is uncertainty in the prior knowledge of $\mtx{A}$, the certificate of
accuracy idea can again be used to authenticate the final answer.

\section{Randomized algorithms for computing a rank-revealing QR decomposition}
\label{sec:randQR}

Up until now, all methods discussed have concerned the problems of computing a
low-rank approximation to a given matrix. These methods were designed explicitly
for the case where the rank $k$ is substantially smaller than the matrix dimensions
$m$ and $n$. In the last two sections, we will describe some recent developments that
illustrate how randomized projections can be used to accelerate matrix factorization
algorithms for any rank $k$, including \textit{full} factorizations where $k=\min(m,n)$.
These new algorithms offer ``one stop shopping'' in that they are faster than
traditional algorithms in essentially every computational regime. We start in this section
with a randomized algorithm for computing full and partial column pivoted QR (CPQR) factorizations.

The material in this section assumes that the reader is familiar with the classical
Householder QR factorization procedure, and with the concept of \textit{blocking}
to accelerate matrix factorization algorithms. It is intended
as a high-level introduction to randomized algorithms for computing full factorizations
of matrices. For details, we refer the reader to \cite{2015_blockQR_SISC,2015_blockQR}.

\subsection{Column pivoted QR decomposition}
\begin{figure}
\setlength{\unitlength}{1mm}
\begin{picture}(125,28)
\put(000,04){\includegraphics[width=24mm]{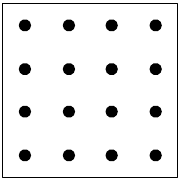}}
\put(032,04){\includegraphics[width=24mm]{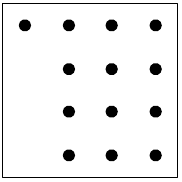}}
\put(064,04){\includegraphics[width=24mm]{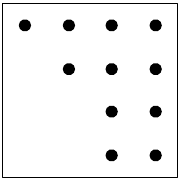}}
\put(096,04){\includegraphics[width=24mm]{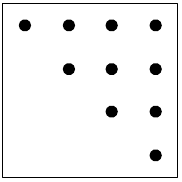}}
\put(007,00){\footnotesize $\mtx{A}_{0} = \mtx{A}$}
\put(035,00){\footnotesize $\mtx{A}_{1} = \mtx{Q}_{1}^{*}\mtx{A}_{0}\mtx{P}_{1}$}
\put(067,00){\footnotesize $\mtx{A}_{2} = \mtx{Q}_{2}^{*}\mtx{A}_{1}\mtx{P}_{2}$}
\put(099,00){\footnotesize $\mtx{A}_{3} = \mtx{Q}_{3}^{*}\mtx{A}_{2}\mtx{P}_{3}$}
\end{picture}
\caption{Basic QR factorization process. The $n\times n$ matrix $\mtx{A}$ is driven to upper
triangular form in $n-1$ steps. (Shown for $n=4$.) At step $i$, we form $\mtx{A}_{i} = \mtx{Q}_{i}^{*}\mtx{A}_{i-1}\mtx{P}_{i}$
where $\mtx{P}_{i}$ is a permutation matrix that moves the largest column in $\mtx{A}_{i-1}(:,1:n)$ to the
$i$'th position, and $\mtx{Q}_{i}$ is a Householder reflector that zeros out all elements under the diagonal
in the pivot column.}
\label{fig:qr}
\end{figure}

Given an $m\times n$ matrix $\mtx{A}$, with $m \geq n$, we recall that the column pivoted QR
factorization (CPQR) takes the form, cf.~(\ref{eq:defCPQS}),
\begin{equation}
\label{eq:defCPQR}
\begin{array}{ccccccccccccccccccc}
\mtx{A}&\mtx{P} &=& \mtx{Q}&\mtx{R},\\
m\times n & n\times n && m\times n & n\times n
\end{array}
\end{equation}
where $\mtx{P}$ is a permutation matrix, where $\mtx{Q}$ is an orthogonal matrix, and where
$\mtx{R}$ is upper triangular. A standard technique for computing the CPQR of a matrix
is the Householder QR process, which we illustrate in Figure
\ref{fig:qr}. The process requires $n-1$ steps to drive $\mtx{A}$ to upper triangular form
from a starting point of $\mtx{A}_{0} = \mtx{A}$. At the $i$'th step, we form
$$
\mtx{A}_{i} = \mtx{Q}_{i}^{*}\mtx{A}_{i-1}\mtx{P}_{i}
$$
where $\mtx{P}_{i}$ is a permutation matrix that flips the $i$'th column of $\mtx{A}_{i-1}$
with the column in $\mtx{A}_{i-1}(:,i:n)$ that has the largest magnitude. The column moved
into the $i$'th place is called the \textit{pivot column.} The matrix $\mtx{Q}_{i}$ is a so
called \textit{Householder reflector} that zeros out all elements beneath the diagonal in
the pivot column. In other words, if we let $\vct{c}_{i}$ denote the pivot column, then
\footnote{The matrix $\mtx{Q}_{i}$ is in fact symmetric, so $\mtx{Q}_{i} = \mtx{Q}_{i}^{*}$, but we
keep the transpose symbol in the formula for consistency with the remainder of the section.}
$$
\mtx{Q}_{i}^{*}\vct{c}_{i} = \left[\begin{array}{c}\vct{c}_{i}(1:(i-1)) \\ \pm \|\vct{c}_{i}(i:m)\| \\ \vct{0}\end{array}\right].
$$
Once the process has completed, the matrices $\mtx{Q}$ and $\mtx{P}$ in (\ref{eq:defCPQR}) are
given by
$$
\mtx{Q} = \mtx{Q}_{n-1}\mtx{Q}_{n-2}\cdots\mtx{Q}_{1},
\qquad\mbox{and}\qquad
\mtx{P} = \mtx{P}_{n-1}\mtx{P}_{n-2}\cdots\mtx{P}_{1}.
$$
For details, see, e.g.,  \cite[Sec.~5.2]{golub}.

\begin{figure}
\setlength{\unitlength}{1mm}
\begin{picture}(125,28)
\put(000,04){\includegraphics[width=24mm]{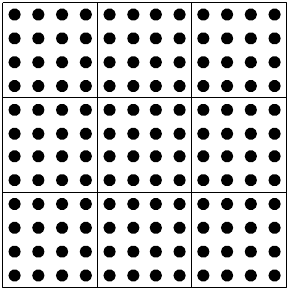}}
\put(032,04){\includegraphics[width=24mm]{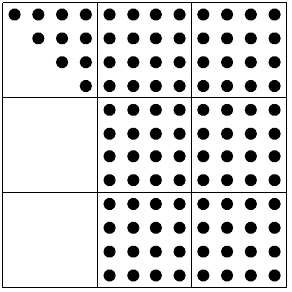}}
\put(064,04){\includegraphics[width=24mm]{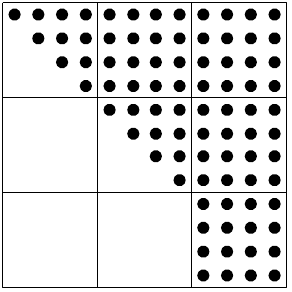}}
\put(096,04){\includegraphics[width=24mm]{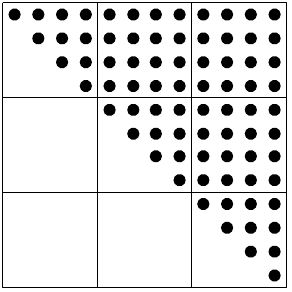}}
\put(007,00){\footnotesize $\mtx{A}_{0} = \mtx{A}$}
\put(035,00){\footnotesize $\mtx{A}_{1} = \mtx{Q}_{1}^{*}\mtx{A}_{0}\mtx{P}_{1}$}
\put(067,00){\footnotesize $\mtx{A}_{2} = \mtx{Q}_{2}^{*}\mtx{A}_{1}\mtx{P}_{2}$}
\put(099,00){\footnotesize $\mtx{A}_{3} = \mtx{Q}_{3}^{*}\mtx{A}_{2}\mtx{P}_{3}$}
\end{picture}
\caption{Blocked QR factorization process. The matrix $\mtx{A}$ consists of $p\times p$ blocks
of size $b\times b$ (shown for $p=3$ and $b=4$). The matrix is driven to upper
triangular form in $p$ steps. At step $i$, we form $\mtx{A}_{i} = \mtx{Q}_{i}^{*}\mtx{A}_{i-1}\mtx{P}_{i}$
where $\mtx{P}_{i}$ is a permutation matrix, and $\mtx{Q}_{i}$ is a product of $b$ Householder reflectors.}
\label{fig:blockqr}
\end{figure}

The Householder QR factorization process is a celebrated algorithm that is exceptionally
stable and accurate. However, it has a serious weakness in that it executes rather slowly
on modern hardware, in particular on systems involving many cores, when the matrix is
stored in distributed memory or on a hard drive, etc. The problem is that it inherently
consists of a sequence of $n-1$ rank-1 updates (so called BLAS2 operations), which makes
the process very communication intensive. In principle, the resolution to this problem
is to \textit{block} the process, as shown in Figure \ref{fig:blockqr}. Let $b$ denote a
block size; then in a blocked Householder QR algorithm, we would find groups of $b$ pivot vectors that are
moved into the active $b$ slots at once, then $b$ Householder reflectors would be
determined by processing the $b$ pivot columns, and then the remainder of the matrix
would be updated jointly. Such a blocked algorithm would expend most of its flops
on matrix-matrix multiplications (so called BLAS3 operations), which execute very
rapidly on a broad range of computing hardware. Many techniques for blocking Householder
QR have been proposed over the years, including, e.g., \cite{1998_quintana_panelRRQR1,1998_quintana_panelRRQR2}.

It was recently observed \cite{2015_blockQR_SISC} that randomized sampling is ideally
suited for resolving the long-standing problem of how to find groups of pivot vectors. The key
observation is that a measure of quality for a group of $b$ pivot vectors is its spanning
volume in $\mathbb{R}^{m}$. This turns out to be closely related to how good of a basis
these vectors form for the column space of the matrix \cite{tyrt1997,gu1996,2006_martinsson_skeletonization}.
As we saw in Section \ref{sec:ID}, this task is particularly well suited to randomized
sampling. To be precise, consider the task of identifying a group of $b$ good pivot vectors
in the first step of the blocked QR process shown in Figure \ref{fig:blockqr}. Using the
procedures described in Section \ref{sec:randID}, we proceed as follows: Fix an over-sampling
parameter $p$, say $p=10$. Then draw a Gaussian random matrix $\mtx{G}$ of size $(b+p)\times m$,
and form the sampling matrix $\mtx{Y} = \mtx{G}\mtx{A}$. Then simply perform column pivoted
QR on the columns of $\mtx{Y}$. To summarize, we determine $\mtx{P}_{1}$ as follows:
\begin{align*}
\mtx{G} =&\ \texttt{randn}(b+p,m),\\
\mtx{Y} =&\ \mtx{G}\mtx{A},\\
[\sim,\sim,\mtx{P}_{1}] =&\ \texttt{qr}(\mtx{Y},0).
\end{align*}
Observe that the QR factorization of $\mtx{Y}$ is affordable since $\mtx{Y}$ is small, and
fits in fast memory close to the processor. For the remaining steps, we simply apply the
same idea to find the best spanning columns for the lower right block in $\mtx{A}_{i}$
that has not yet been driven to upper triangular form. The resulting algorithm is called
\textit{Householder QR with Randomization for Pivoting (HQRRP)}; it is described in detail
in \cite{2015_blockQR_SISC}, and is available at \texttt{https://github.com/flame/hqrrp/}.
(The method described in this section was first published in \cite{2015_blockQR}, but is closely
related to the independently discovered results in \cite{2015_blockQR_ming}.)

To maximize performance, it turns out to be possible to ``downdate'' the
sampling matrix from one step of the factorization to the next, in a manner
similar to how downdating of the pivot weights are done in classical
Householder QR\cite[Ch.5, Sec.~2.1]{1998_stewart_volume1}. This obviates
the need to draw a new random matrix at each step \cite{2015_blockQR_ming,2015_blockQR_SISC},
and reduces the leading term in the asymptotic flop count of HQRRP to $2mn^{2} - (4/3)n^{3}$, which
is identical to classical Householder QR.

\section{A strongly rank-revealing UTV decomposition}
\label{sec:randUTV}


This section describes a randomized algorithm \texttt{randUTV} that is very similar to the randomized
QR factorization process described in Section \ref{sec:randQR} but that results in a so called
``UTV factorization.'' The new algorithm has several advantages:
\begin{itemize}
\item \texttt{randUTV} provides close to optimal low-rank approximation, and highly accurate
estimates for the singular values of a matrix.
\item The algorithm \texttt{randUTV} builds the factorization (\ref{eq:defUTVpre})
\textit{incrementally,}
which means that when it is applied to a matrix of numerical rank $k$, the algorithm can be
stopped early and incur an overall cost of $O(mnk)$.
\item Like \texttt{HQRRP}, the algorithm \texttt{randUTV} is \textit{blocked}, which enables
it to execute fast on modern hardware.
\item The algorithm \texttt{randUTV} is not an iterative algorithm. In this regard, it is closer
to the CPQR than standard SVD algorithms, which substantially simplifies software optimization.
\end{itemize}

\subsection{The UTV decomposition}
Given an $m\times n$ matrix $\mtx{A}$, with $m\geq n$, a ``UTV decomposition'' of $\mtx{A}$ is a
factorization the form
\begin{equation}
\label{eq:defUTVpre}
\begin{array}{ccccccccccc}
\mtx{A} &=& \mtx{U} & \mtx{T} & \mtx{V}^{*},\\
m\times n && m\times m & m\times n & n\times n
\end{array}
\end{equation}
where $\mtx{U}$ and $\mtx{V}$ are unitary matrices, and $\mtx{T}$ is a triangular matrix
(either lower or upper triangular).
The UTV decomposition can be viewed as a generalization of other standard factorizations
such as, e.g., the \textit{Singular Value Decomposition (SVD)} or the
\textit{Column Pivoted QR decomposition (CPQR).} (To be precise, the SVD is the special case
where $\mtx{T}$ is diagonal, and the CPQR is the special case where $\mtx{V}$ is a permutation matrix.)
The additional flexibility inherent in the UTV
decomposition enables the design of efficient updating procedures,
see \cite[Ch.~5, Sec.~4]{1998_stewart_volume1} and \cite{1999_hansen_UTVtools}.

\subsection{An overview of \texttt{randUTV}}
\label{sec:UTVoverview}
\begin{figure}
\setlength{\unitlength}{1mm}
\begin{picture}(125,28)
\put(000,04){\includegraphics[width=24mm]{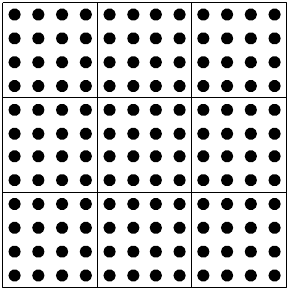}}
\put(032,04){\includegraphics[width=24mm]{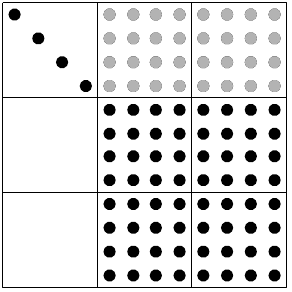}}
\put(064,04){\includegraphics[width=24mm]{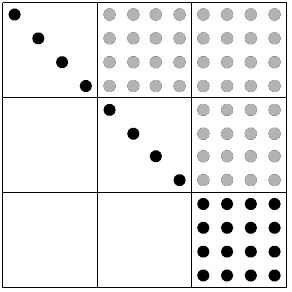}}
\put(096,04){\includegraphics[width=24mm]{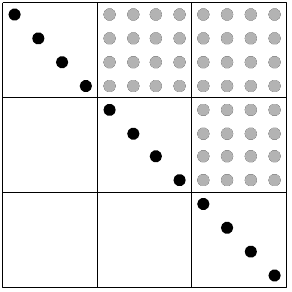}}
\put(007,00){\footnotesize $\mtx{A}_{0} = \mtx{A}$}
\put(035,00){\footnotesize $\mtx{A}_{1} = \mtx{U}_{1}^{*}\mtx{A}_{0}\mtx{V}_{1}$}
\put(067,00){\footnotesize $\mtx{A}_{2} = \mtx{U}_{2}^{*}\mtx{A}_{1}\mtx{V}_{2}$}
\put(099,00){\footnotesize $\mtx{A}_{3} = \mtx{U}_{3}^{*}\mtx{A}_{2}\mtx{V}_{3}$}
\end{picture}
\caption{Blocked UTV factorization process. The matrix $\mtx{A}$ consists of $p\times p$ blocks
of size $b\times b$ (shown for $p=3$ and $b=4$). The matrix is driven to upper
triangular form in $p$ steps. At step $i$, we form $\mtx{A}_{i} = \mtx{U}_{i}^{*}\mtx{A}_{i-1}\mtx{V}_{i}$
where $\mtx{U}_{i}$ and $\mtx{V}_{i}$ consist (mostly) of a product of $b$ Householder reflectors.
The elements shown in grey are not zero, but are very small in magnitude.}
\label{fig:blockUTV}
\end{figure}

The algorithm \texttt{randUTV} follows the same general pattern as \texttt{HQRRP}, as illustrated
in Figure \ref{fig:blockUTV}. Like \texttt{HQRRP}, it drives a given matrix $\mtx{A} = \mtx{A}_{0}$
to upper triangular form via a sequence of steps
$$
\mtx{A}_{i} = \mtx{U}_{i}^{*}\mtx{A}_{i-1}\mtx{V}_{i},
$$
where each $\mtx{U}_{i}$ and $\mtx{V}_{i}$ is a unitary matrix.
As in Section \ref{sec:randQR}, we let $b$ denote a block size, and let $p = \lceil n/b\rceil$
denote the number of steps taken.
The key difference from \texttt{HQRRP} is that we in \texttt{randUTV} allow the matrices
$\mtx{V}_{i}$ to consist (mostly) of a  product of $b$ Householder reflectors. This added
flexibility over CPQR allows us to drive more mass onto the diagonal entries, and thereby
render the off-diagonal entries in the final matrix $\mtx{T}$ very small in magnitude.

\subsection{A single step block factorization}
\label{sec:stepUTV}
The algorithm \texttt{randUTV} consists of repeated application of a randomized
technique for building approximations to the spaces spanned by the domaninant
$b$ left and right singular vectors, where $b$ is a given block size. To be precise,
given an $m\times n$ matrix $\mtx{A}$, we seek to build unitary matrices $\mtx{U}_{1}$
and $\mtx{V}_{1}$ such that
$$
\mtx{A} = \mtx{U}_{1}\mtx{A}_{1}\mtx{V}_{1}^{*}
$$
where $\mtx{A}_{1}$ has the block structure
$$
\mtx{A}_{1} =
\left[\begin{array}{cc}
\mtx{A}_{1,11} & \mtx{A}_{1,12} \\
\mtx{0}        & \mtx{A}_{1,22}
\end{array}\right] =
\raisebox{-9mm}{\includegraphics[height=18mm]{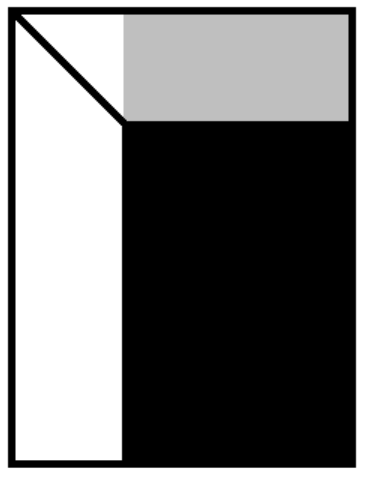}},
$$
so that $\mtx{A}_{1,11}$ is diagonal, and $\mtx{A}_{1,12}$ has entries
of small magnitude.

We first build $\mtx{V}_{1}$. To this end, we use the randomized power iteration
described in Section \ref{sec:power} to build a sample matrix $\mtx{Y}$ of size
$b\times n$ whose columns approximately span the same subspace as the $b$
dominant right singular vectors of $\mtx{A}$. To be precise, we draw a $b\times m$ Gaussian
random matrix $\mtx{G}$ and form the sample matrix
$$
\mtx{Y} = \mtx{G}\mtx{A}\bigl(\mtx{A}^{*}\mtx{A}\bigr)^{q},
$$
where $q$ is a parameter indicating the number of steps of power iteration taken
(in \texttt{randUTV}, the gain from over-sampling is minimal and is generally not
worth the bother). Then we form a unitary
matrix $\tilde{\mtx{V}}$ whose first $b$ columns form an orthonormal basis for the
column space of $\mtx{Y}$. (The matrix $\mtx{V}$ consists of a product of $b$ Householder
reflectors, which are determined by executing the standard Householder QR procedure on
the columns of $\mtx{Y}^{*}$.) We then execute $b$ steps of Householder QR
on the matrix $\mtx{A}\tilde{\mtx{V}}$ to form a matrix $\tilde{\mtx{U}}$ consisting of
a product of $b$ Householder reflectors. This leaves us with a new matrix
$$
\tilde{\mtx{A}} = \bigl(\tilde{\mtx{U}}\bigr)^{*}\,\mtx{A}\,\tilde{\mtx{V}}
$$
that has the block structure
$$
\tilde{\mtx{A}} =
\left[\begin{array}{cc}
\tilde{\mtx{A}}_{11} & \tilde{\mtx{A}}_{12} \\
\mtx{0}        & \tilde{\mtx{A}}_{22}
\end{array}\right] =
\raisebox{-9mm}{\includegraphics[height=18mm]{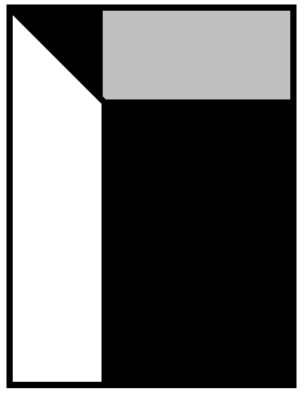}}.
$$
In other words, the top left $b\times b$ block is upper triangular, and the bottom left
block is zero. One can also show that all entries of $\tilde{\mtx{A}}_{12}$ are typically
small in magnitude.
Next, we compute a full SVD of the block $\tilde{\mtx{A}}_{11}$
$$
\tilde{\mtx{A}}_{11} = \hat{\mtx{U}}\mtx{D}_{11}\hat{\mtx{V}}^{*}.
$$
This step is affordable since $\tilde{\mtx{A}}_{11}$ is of size $b\times b$, where $b$ is small.
Then form the transformation matrices
$$
\mtx{U}_{1} =
\tilde{\mtx{U}}
\left[\begin{array}{cc}
\hat{\mtx{U}} & \mtx{0} \\
\mtx{0} & \mtx{I}_{m-b}
\end{array}\right],
\qquad\mbox{and}\qquad
\mtx{V}_{1} =
\tilde{\mtx{V}}
\left[\begin{array}{cc}
\hat{\mtx{V}} & \mtx{0} \\
\mtx{0} & \mtx{I}_{n-b}
\end{array}\right],
$$
and set
$$
\mtx{A}_{1} = \mtx{U}_{1}^{*}\mtx{A}\mtx{V}_{1}.
$$
One can demonstrate that the diagonal entries of $\mtx{D}_{11}$ typically
form  accurate approximations to first $b$ singular values of $\mtx{A}$, and
that
$$
\|\mtx{A}_{1,22}\| \approx \inf\{\|\mtx{A} - \mtx{B}\|\,\colon\,\mtx{B}\mbox{ has rank }b\}.
$$

Once the first $b$ columns and rows of $\mtx{A}$ have been processed as described in this
Section, \texttt{randUTV} then applies the same procedure
to the remaining block $\mtx{A}_{1,22}$, which has size $(m-b) \times (n-b)$,
and then continues in the same fashion to process all remaining blocks,
as outlined in Section \ref{sec:UTVoverview}. The full algorithm is summarized in Figure \ref{fig:randUTV}.

For more information about the UTV factorization, including careful numerical experiments
that illustrate how it compares in terms of speed and accuracy to competitors such as
column pivoted QR and the traditional SVD, see \cite{2017_martinsson_UTV}. Codes are
available at \texttt{https://github.com/flame/randutv}.

\begin{remark}
As discussed in Section \ref{sec:power}, the rows of $\mtx{Y}$ approximately span
the linear space spanned by the $b$ dominant right singular vectors of $\mtx{A}$.
The first $b$ columns of $\mtx{V}_{1}$ simply form an orthonormal basis for
$\mbox{Row}(\mtx{Y}^{*})$. Analogously, the first $b$ columns of $\mtx{U}_{1}$
approximately form an orthonormal basis for the $b$ dominant left singular vectors
of $\mtx{A}$. However, the additional application of $\mtx{A}$ that is implicit in
forming the matrix $\mtx{A}\tilde{\mtx{V}}$ provides a boost in accuracy, and the
rank-$b$ approximation resulting from one step of \texttt{randUTV} is similar to the
accuracy obtained from the randomized algorithm in Section \ref{sec:power}, but with
``$q+1/2$ steps'' of power iteration, rather than $q$ steps.
(We observe that if the first $b$ columns of $\tilde{\mtx{U}}$ and $\tilde{\mtx{V}}$
spanned \textit{exactly} the spaces spanned by the dominant $b$ left and right singular
vectors of $\mtx{A}$, then we would have $\tilde{\mtx{A}}_{12} = \mtx{0}$, and
the singular values of $\tilde{\mtx{A}}_{11}$ would be identical to the top $b$ singular
values of $\mtx{A}$.)
\end{remark}

\begin{figure}
\begin{center}
\fbox{\begin{minipage}{170mm}
\begin{tabbing}
\hspace{8mm} \= \hspace{5mm} \= \hspace{5mm} \= \hspace{15mm}\kill
\textbf{function} $[\mtx{U},\mtx{T},\mtx{V}] = \texttt{randUTV}(A,b,q)$\\
\> $\mtx{T} = \mtx{A}$;\\
\> $\mtx{U} = \texttt{eye}(\texttt{size}(\mtx{A},1))$;\\
\> $\mtx{V} = \texttt{eye}(\texttt{size}(\mtx{A},2))$;\\
\> \textbf{for} $i = 1:\texttt{ceil}(\texttt{size}(\mtx{A},2)/b)$\\
\> \> $I_{1} = 1:(b(i-1))$;\\
\> \> $I_{2} = (b(i-1)+1):\texttt{size}(\mtx{A},1)$;\\
\> \> $J_{2} = (b(i-1)+1):\texttt{size}(\mtx{A},2)$;\\
\> \> \textbf{if} $(\texttt{length}(J_{2}) > b)$\\
\> \> \> $[\hat{\mtx{U}},\hat{\mtx{T}},\hat{\mtx{V}}] = \texttt{stepUTV}(\mtx{T}(I_{2},J_{2}),b,q)$;\\
\> \> \textbf{else}\\
\> \> \> $[\hat{\mtx{U}},\hat{\mtx{T}},\hat{\mtx{V}}] = \texttt{svd}(\mtx{T}(I_{2},J_{2}))$;\\
\> \> \textbf{end if}\\
\> \> $\mtx{U}(:,I_{2})  = \mtx{U}(:,I_{2})*\hat{\mtx{U}}$;\\
\> \> $\mtx{V}(:,J_{2})  = \mtx{V}(:,J_{2})*\hat{\mtx{V}}$;\\
\> \> $\mtx{T}(I_{2},J_{2}) = \hat{\mtx{T}}$;\\
\> \> $\mtx{T}(I_{1},J_{2}) = \hat{T}(I_{1},J_{2})*\hat{\mtx{V}}$;\\
\> \textbf{end for}\\
\texttt{return}
\end{tabbing}
\end{minipage}}

\fbox{\begin{minipage}{170mm}
\begin{tabbing}
\hspace{8mm} \= \hspace{5mm} \= \hspace{5mm} \= \hspace{15mm}\kill
\textbf{function} $[U,T,V] = \texttt{stepUTV}(A,b,q)$\\
\> $\mtx{G} = \texttt{randn}(\texttt{size}(A,1),b)$;\\
\> $\mtx{Y} = \mtx{A}^{*}\mtx{G}$;\\
\> \textbf{for} $i = 1:q$\\
\>\> $\mtx{Y} = \mtx{A}^{*}(\mtx{A}\mtx{Y})$;\\
\> \textbf{end for}\\
\> $[\mathbb{V},\sim] = \texttt{qr}(Y)$;\\
\> $[\mtx{U},\mtx{D},\mtx{W}] = \texttt{svd}(\mtx{A}\mtx{V}(:,1:b))$;\\
\> $\mtx{T} = \bigl[\mtx{D},\ \mtx{U}^{*}\mtx{A}\mtx{V}(:,(b+1):\texttt{end})\bigr]$;\\
\> $\mtx{V}(:,1:b) = \mtx{V}(:,1:b)*\mtx{W}$;\\
\textbf{return}
\end{tabbing}
\end{minipage}}
\end{center}
\caption{The algorithm \texttt{randUTV} (described in Section \ref{sec:randUTV}) that given an $m\times n$ matrix $\mtx{A}$
computes its UTV factorization $\mtx{A} = \mtx{U}\mtx{T}\mtx{V}^{*}$, cf.~(\ref{eq:defUTVpre}).
The input parameters $b$ and $q$ reflect the block size and the number of steps of power iteration,
respectively. The single step function \texttt{stepUTV} is described in Section \ref{sec:stepUTV}.
(Observe that most of the unitary matrices that arise consist of products of Householder reflectors;
this property must be exploited to attain computational efficiency.)}
\label{fig:randUTV}
\end{figure}

\clearpage

\bibliographystyle{plain}
\bibliography{main_bib}

\end{document}